# The Pareto-Optimal Temporal Aggregation of Energy System Models

Maximilian Hoffmann,[a,b,1] Leander Kotzur,[a] and Detlef Stolten[a,b]

[a] Institute of Energy and Climate Research, Techno-economic Systems Analysis (IEK-3), Forschungszentrum Jülich, 52428 Jülich, Germany

[b] Chair for Fuel Cells, RWTH Aachen University, c/o Institute of Electrochemical Process Engineering (IEK-3), Forschungszentrum Jülich GmbH, Wilhelm-Johnen-Str., 52428 Jülich, Germany

## Abstract

The growing share of intermittent renewable energy sources, storage technologies, and the increasing degree of so-called sector coupling necessitates optimization-based energy system models with high temporal and spatial resolutions, which significantly increases their runtimes and limits their maximum sizes. In order to maintain the computational viability of these models for large-scale application cases, temporal aggregation has emerged as a technique for reducing the number of considered time steps by reducing the original time horizon down to fewer, more representative ones.

This study presents advanced but generally applicable clustering techniques that allow for ad-hoc improvements of state-of-the-art approaches without requiring profound knowledge of the individual energy system model. These improvements comprise the optimal tradeoff between the number of typical days and inner-daily temporal resolutions, as well as constituting a representation method that can reproduce the value distribution of the original time series. We prove the superiority of these approaches by applying them to two fundamentally different model types, namely a single-node building energy system and a European carbon-neutral energy scenario, and benchmark these against state-of-the-art approaches. This is performed for a variety of temporal resolutions, which leads to many hundreds of model runs.

The results show that the proposed improvements on current methods strictly dominate the status quo with respect to Pareto-optimality in terms of runtime and accuracy. Although a speeding up factor of one magnitude could be achieved using traditional aggregation methods within a cost deviation range of two percent, the algorithms proposed herein achieve this accuracy with a runtime speedup by a factor of two orders of magnitude.

## Highlights

- The ratio of typical days and time steps therein is crucial for a good optimization result.
- An algorithm for an improved ratio of typical days and segments is introduced.
- The duration curves of time series are crucial for the correct sizing of the components.
- An algorithm is introduced that precisely replicates the original value distribution.
- The methods are open source and available at: https://github.com/FZJ-IEK3-VSA/tsam.

[1] Corresponding author: M. Hoffmann (max.hoffmann@fz-juelich.de)

# 1. Introduction

Current and future energy systems will be subject to a steadily increasing share of intermittent renewable energy sources due to the necessity of reducing carbon dioxide emissions, which also leads to an intensification of synergies, or 'coupling', across different energy sectors [1, 2]. These trends increase system complexity and drive the need for spatially- and temporally-highly-resolved energy system models. In contrast, the growth rate of transistor density, which has been constant over the decades and is known as Moore's Law [3, 4], has been decreasing in recent years [5].

The conflict between fast-growing energy system models and a decreasing growth rate of computing power has driven the development of algorithms that focus on reducing the complexity of energy system models [6, 7]. One of these approaches has been the reduction of the number of time steps considered in these time-discrete models, for which a large variety of different approaches exist [8]. Among these are aggregation procedures based on clustering, which have emerged as the most accurate methods due to their capability of accounting for the similarity of different time steps or periods. However, most methods described in the literature either focus on reducing the number of time steps by increasing their duration or by representing the time series with a few typical periods, whereas the use of both methods together has been widely neglected. Furthermore, the vast majority of studies have employed well-established clustering and representation methods, such as k-means, k-medoids, and Ward's hierarchical algorithm, all of which have led to a distortion of the original time series' extreme values and value distributions (duration curve). Accordingly, this work addresses two major research questions:

**Research question 1:** Is there an optimal tradeoff between the duration of time steps and the number of typical periods, and are certain technologies affected systematically by either of these temporal aggregation approaches?

**Research question 2:** Can the aggregation procedure be improved upon by explicitly taking statistical momentums, i.e., the original time series' value distribution, into account?

Both research questions are investigated using two fundamentally different energy system models, and the methods developed in this study are shown to outperform state-of-the-art approaches using the criterion of Pareto-optimality.

The remainder of this work is structured as follows: In Section 2, we review current developments in temporal aggregation for energy system models. Section 3 introduces the modifications of standard aggregation procedures as proposed herein and provides a short overview of the case studies used to validate the effectiveness of the proposed modifications. Section 4 analyzes the impact of the proposed aggregation procedures on the given case study. Finally, the findings are summarized and concluded in Section 5.

# 2. State-of-the-art

The following section is divided into two parts. Section 1.1 reviews approaches to reducing the number of time steps, whereas Section 0 outlines subsequent methods for improving the quality of the aggregated time series by manually or automatically inserting the statistical features of the fully resolved time series that could be relevant for energy system optimization.

## 1.1. Typical periods and segments

Time series can be aggregated in two different ways, either by directly reducing the temporal resolution, i.e., merging adjacent time steps, or by aggregating multiple periods of consecutive time steps into fewer periods. As many time series in the realm of energy system modeling, such as of electricity or photovoltaic feed-in profiles, have a daily pattern, the most common period length is 24 hours. In that case, this method is also referred to as *typical* or *representative days*. As this approach affects the number of considered periods but not the actual temporal resolution of discrete time steps within these, the temporal resolution can be further reduced in a subsequent step. The impact of aggregation on a typical single time series is depicted in Figure 1. Here, it can be seen that the aggregation to typical periods leads to fewer periods, whereas segmentation leads directly to fewer time steps using a coarser resolution. In the example, the original time series with 8760 hourly time steps was reduced to eight typical days comprising eight segments with irregular lengths each, ending up in 64 time steps, which is equal to a time step reduction of well above 99%. As the runtime of energy system models scales nonlinearly with the number of simultaneously considered time steps, this reduction may lead to an even higher speed-up factor.

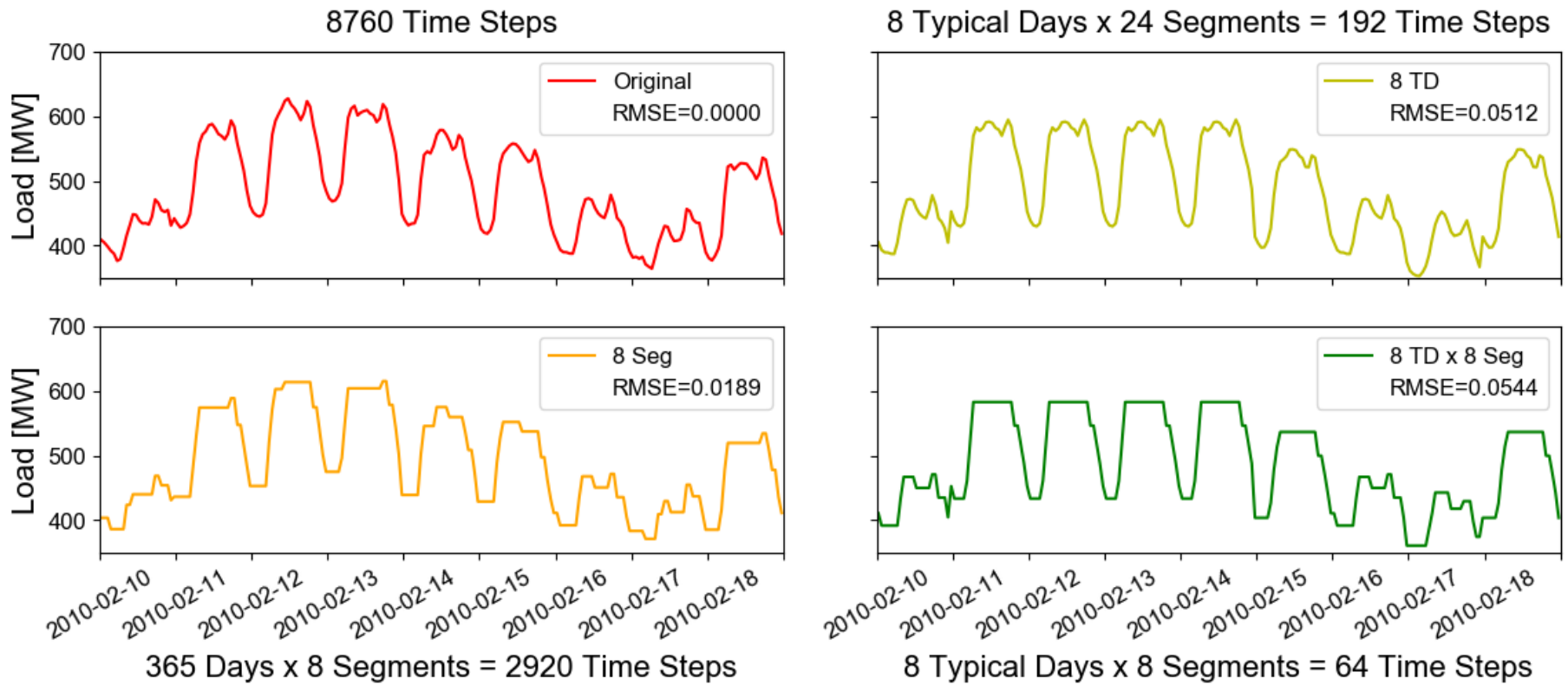

*Figure 1. Time step reduction based on a combination of period clustering and segmentation.*

Despite the option to freely combine typical periods and a decrease in inner-periodical temporal resolution, very few authors have considered both approaches simultaneously in the literature, which is illustrated in Table 1.

The methods that resemble the approach taken in this study presented in Section 3.2 were mostly developed by Fazlollahi et al. [9], as well as Bahl et al. [10-13]. Fazlollahi et al. employed k-means clustering in order to cluster both periods and the time steps within these periods until a set of input data-based quality criteria were fulfilled. As the segmentation must fulfill the constraint of adjacency, which cannot be directly respected by the k-means algorithm, the segmentation is performed in an iterative manner. In contrast, Bahl et al. perform optimization runs for a certain number of typical days and segments per day in order to obtain upper and lower bounds of the original problem. Then, they either increased the number of segments within each typical day or the number of typical days based on a heuristic to select the more promising direction of refinement in order to obtain tighter bounds on the original model. Finally, Raventos et al. [14] compared the approaches of typical days and segments with one another in order to evaluate their advantages and disadvantages, but did not combine both approaches.

The fact that these approaches did not explicitly compare the performance of various cross-combined aggregation configurations drove the motivation for a large-scale sensitivity analysis, as performed in the following.

*Table 1. Studies using either typical periods or segments based on clustering.*

| Publication | Year | Typical Periods | Segments | Both |
|---|---|---|---|---|
| Domínguez-Muñoz et al. [15] | 2011 | ✖ | | |
| Adhau et al. [16] | 2014 | ✖ | | |
| Green et al. [17] | 2014 | ✖ | | |
| Wogrin et al. [18, 19] | 2014/16 | ✖ | | |
| Agapoff et al. [20] | 2015 | ✖ | | |
| Pfenninger et al. [21] | 2017 | ✖ | | |
| Brodrick et al. [22, 23] | 2015/17 | ✖ | | |
| Fitiwi et al. [24] | 2015 | ✖ | | |
| Marquant et al. [25, 26] | 2015/17 | ✖ | | |
| Munoz et al. [27] | 2015 | ✖ | | |
| Schiefelbein et al. [28] | 2015 | ✖ | | |
| Lin et al. [29] | 2016 | ✖ | | |
| Ploussard et al. [30] | 2016 | ✖ | | |
| Nahmmacher et al. [31] | 2016 | ✖ | | |
| Bahl et al. [32] | 2017 | ✖ | | |
| Härtel et al. [33] | 2017 | ✖ | | |
| Heuberger et al. [34] | 2017 | ✖ | | |
| Zhu et al. [35] | 2017 | ✖ | | |
| Stadler et al. [36] | 2018 | ✖ | | |
| Kotzur et al. [37, 38] | 2018 | ✖ | | |
| Lara et al. [39] | 2018 | ✖ | | |
| Liu et al. [40] | 2018 | ✖ | | |
| Tejada-Arango et al. [41, 42] | 2018 | ✖ | | |
| Tupper et al. [43] | 2018 | ✖ | | |
| Welder et al. [44] | 2018 | ✖ | | |
| Schütz et al. [45-47] | 2016-18 | ✖ | | |
| Gabrielli et al. [48, 49] | 2018/19 | ✖ | | |
| Teichgräber et al. [50-53] | 2017-20 | ✖ | | |
| Kannengießer et al. [54] | 2019 | ✖ | | |
| Sun et al. [55] | 2019 | ✖ | | |
| Zatti et al. [56] | 2019 | ✖ | | |
| Zhang et al. [57] | 2019 | ✖ | | |
| Pinel et al. [58] | 2020 | ✖ | | |
| Raventos et al. [14] | 2020 | ✖ | | |
| Hoffmann et al. [59] | 2021 | ✖ | | |
| Teichgräber et al. [60] | 2021 | ✖ | | |
| Raventos et al. [14] | 2020 | | ✖ | |
| Deml et al. [61] | 2015 | | ✖ | |
| Vom Stein et al. [62] | 2017 | | ✖ | |
| Pineda et al. [63] | 2018 | | ✖ | |
| Mavrotas et al. [64] | 2008 | ✖ | ✖ | ✖ |
| Fazlollahi et al. [9] | 2014 | ✖ | ✖ | ✖ |
| Bahl et al. [10, 11] | 2018 | ✖ | ✖ | ✖ |
| Baumgärtner et al. [12, 13] | 2019 | ✖ | ✖ | ✖ |

## 1.2. Modified clustering techniques

Apart from clustering, time series aggregation also employs a representation step as highlighted in Figure 2. In general, arbitrary value tuples can be used to represent a cluster of time steps from the original time series. However, medoids or centroids are the most frequently used in the literature.

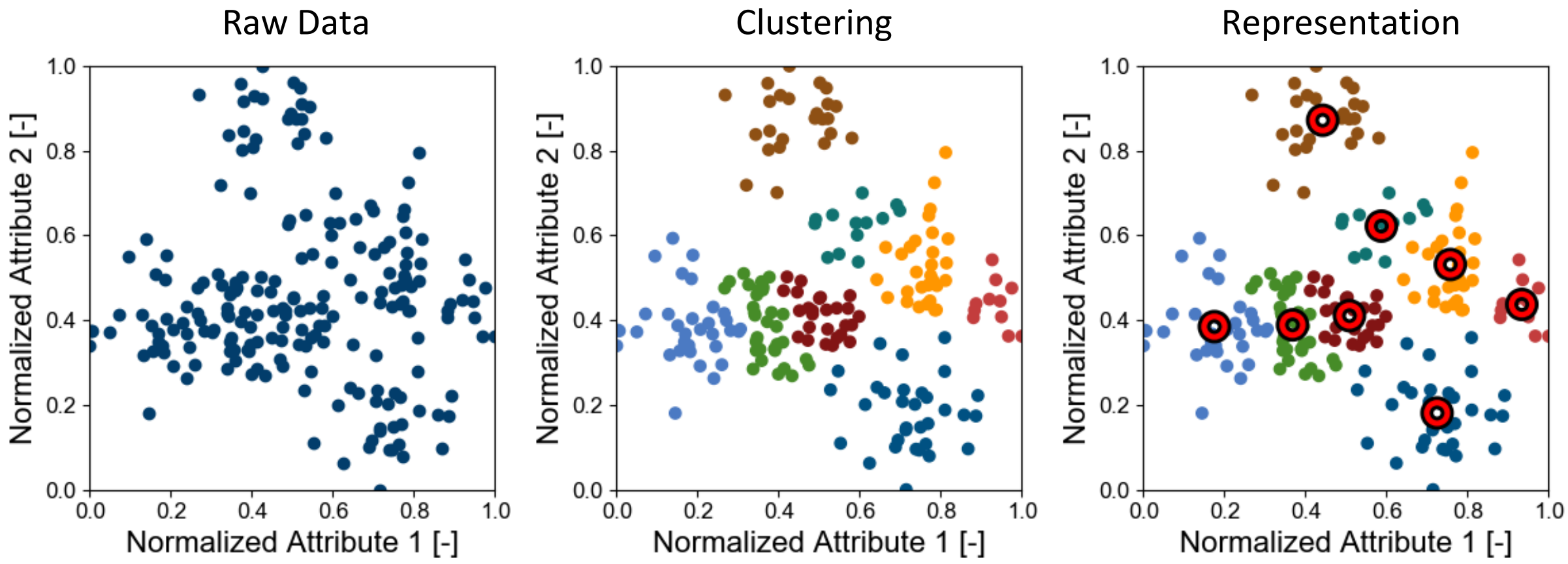

*Figure 2. Clustering and representation as separate steps of clustering.*

The problem of using centroids or medoids is illustrated in Figure 3 for the case of centroids. It can be seen that centroids lead to an underestimation of the time series' extreme values, which is apparent for both its duration curve (left) and an extraction of the original time series during a week in February. This underestimation often leads to an underestimation of the system costs as demand maxima are underestimated and supply minima overestimated.

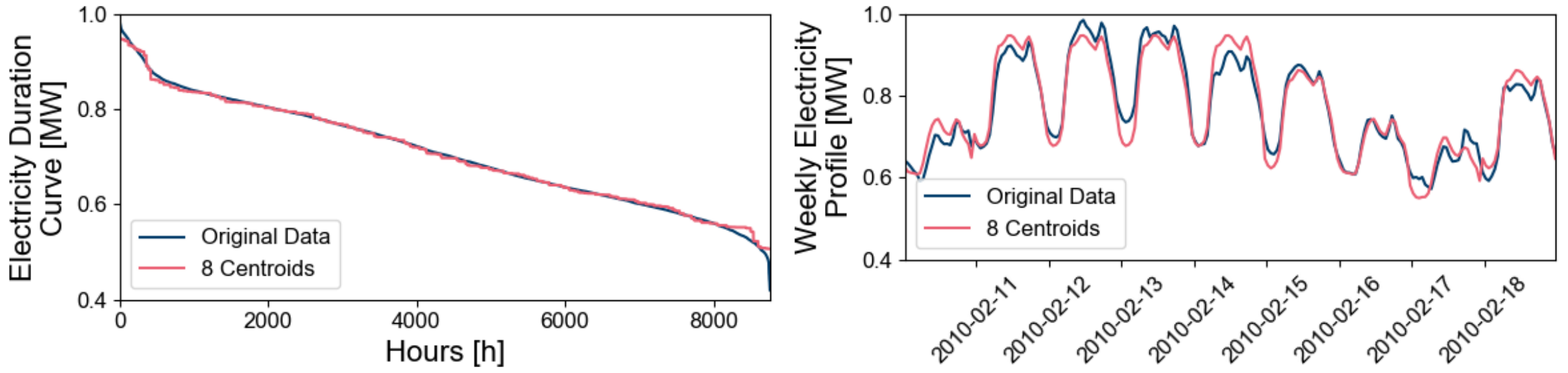

*Figure 3. The yearly duration curve of a typical electricity profile, its profile during a week in February, and the corresponding aggregated time series as predicted by eight typical days using hierarchical clustering and centroids as representative values.*

To address this issue, many approaches have been developed in the literature that either are based on the manual integration of extreme values, or are algorithms that automatically preserve the desirable features of the original time series. Table 2 provides an overview of the different approaches.

*Table 2. Methods employed to subsequently improve the quality of the clustering process.*

| Publication | Year | Manual Peak Integration | Algorithm |
|---|---|---|---|
| Mavrotas et al. [64] | 2008 | ✖ | |
| Domínguez-Muñoz et al. [15] | 2011 | ✖ | |
| Ortiga et al. [65] | 2011 | ✖ | |
| Devogelaer et al. [66] | 2012 | ✖ | |
| De Sisternes et al. [67, 68] | 2013/16 | ✖ | |
| Simões et al. [69] | 2013 | ✖ | |
| Voll et al. [70] | 2013 | ✖ | |
| Fazlollahi et al. [9, 71] | 2014 | ✖ | |
| Green et al. [17] | 2014 | | ✖ |
| Poncelet et al. [72, 73] | 2014/16 | ✖ | |
| Stadler et al. [74] | 2014 | ✖ | |
| Wakui et al. [75-77] | 2014/16 | ✖ | |
| Agapoff et al. [20] | 2015 | | ✖ |
| Marquant et al. [25, 26] | 2015 | ✖ | |
| Munoz et al. [27] | 2015 | ✖ | |
| Frew et al. [78] | 2016 | ✖ | |
| Merrick et al. [79] | 2016 | ✖ | |
| Patteeuw et al. [80] | 2016 | ✖ | |
| Poncelet et al. [81] | 2016 | | ✖ |
| Bahl et al. [10, 11, 32] | 2017/18 | | ✖ |
| Härtel et al. [33] | 2017 | | ✖ |
| Heuberger et al. [34] | 2017 | ✖ | |
| Pfenninger et al. [21] | 2017 | ✖ | |
| Brodrick et al. [23, 82] | 2017/18 | ✖ | |
| Gabrielli et al. [48, 49] | 2018/19 | ✖ | |
| Kotzur et al. [38] | 2018 | ✖ | |
| Baumgärtner et al. [12, 13] | 2019 | ✖ | |
| Hilbers et al. [83] | 2019 | | ✖ |
| Zatti et al. [56] | 2019 | | ✖ |
| Teichgräber et al. [60] | 2021 | | ✖ |

Among the algorithm-based approaches, the findings of Agapoff et al. [20] and Poncelet et al. [81] have influenced the method developed in Section 3.3 the most due to their focus on the statistical features of the original time series, rather than single extreme values. Both authors refrained from directly clustering the original time series. Instead, Agapoff et al. selected local price differences, as well as non-controllable demand and generation time series for transmission expansion planning, as well as their statistical features of average, minimum, maximum and standard deviations to directly perform the clustering of statistical attributes instead of the original time series. In contrast, Poncelet et al. introduced a MILP approach for clustering time series in such a way that their respective duration curves are as closely approximated as possible. Based on these approaches, the method presented in Section 3.3 was developed to explicitly preserve the original value distribution and thereby the statistical key features of the original time series.

In summary, current methods offer very little information on a good ratio between the number of typical days and their inner-daily resolutions, or more generally approaches to account for extreme situations, without manually manipulating certain time series on the basis of the modeler's experience.

# 3. Methodology

In the following, the approaches proposed to improve the aggregation process are introduced, as well as the case studies used for their validation.

## 3.1. Clustering typical days and segments

The clustering process used to aggregate the time series to typical days utilized in this work is based on a procedure that was first introduced by Domínguez-Muñoz et al. [15], and which is frequently described in the literature, e.g., by Nahmmacher et al. [31], as well as the authors of this work [8, 38, 59], who also made this functionality available open-source as part of their tsam modeling tool.[2] For the sake of completeness, however, the process is also depicted in Figure 2. In a first step, the time series are normalized in order to prevent an overweighing of those with large scales (e.g., total energy demands) compared to those with small scales (e.g., capacity factors). This can either happen using a min–max normalization to values of between 0 and 1, or by using a z-normalization, which leads to

[2] tsam package available at: https://github.com/FZJ-IEK3-VSA/tsam (date accessed: 11/23/2021)

time series with a standard deviation of 1 [59]. Then, the time series are rearranged such that all steps of all time series within a certain period, e.g., one day, become a row vector each of which can then be interpreted as a sample point in a hyper-dimensional space and therefore clustered with other sample points in order to achieve groups of days that share a certain degree of similarity. In our study, we employed hierarchical clustering [84] due to its deterministic behavior. Ultimately, one representative is determined for each cluster using either the cluster's centroid, medoid, or a user-specified approach and scaled back. In the energy system model, these typical periods are then weighted according to their respective cluster sizes in order to account for the relative frequency of certain day types throughout the considered time horizon.

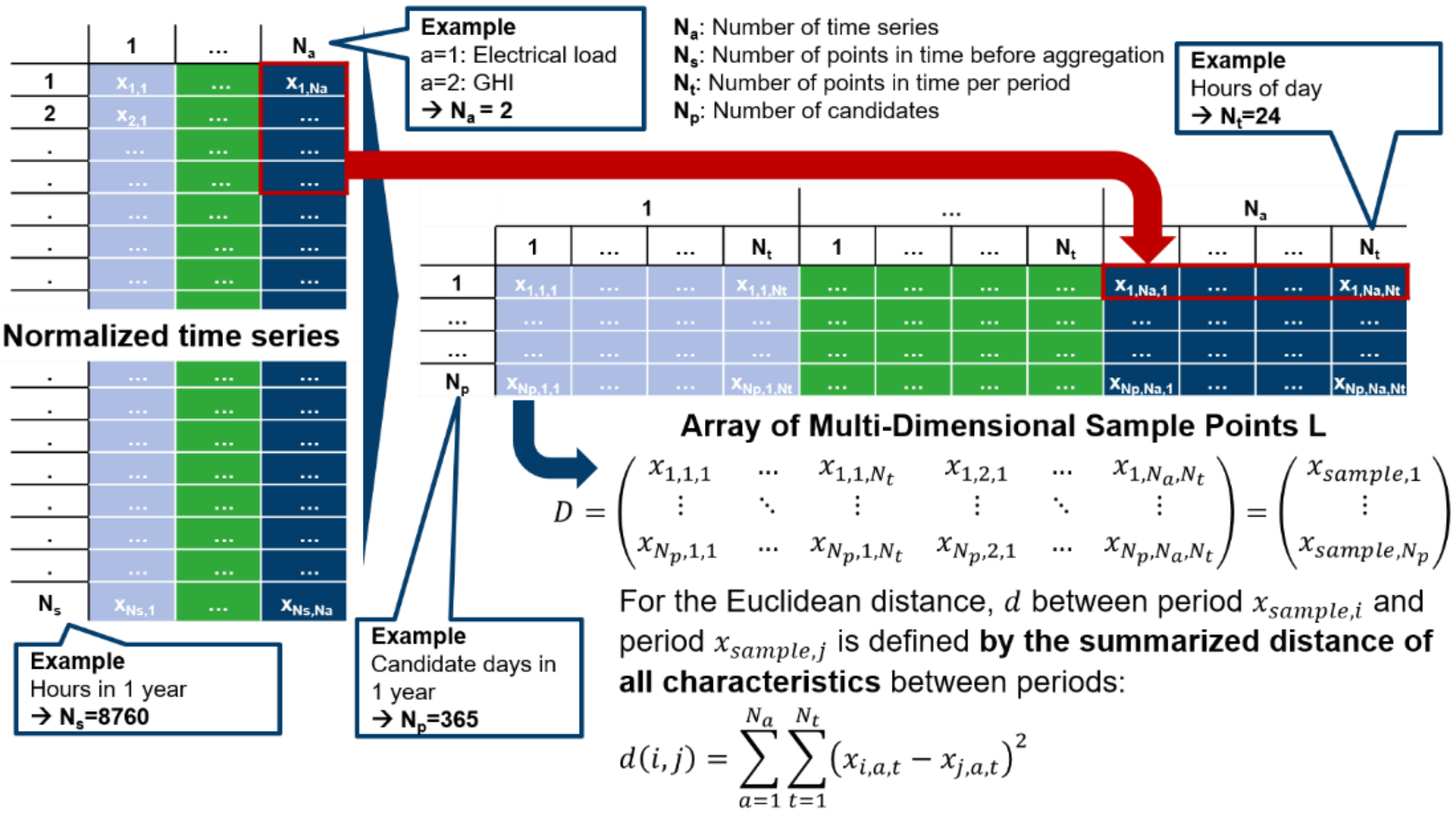

Figure 4. The procedure of period clustering using an arbitrary clustering algorithm.

The method of segmentation applied in this work is also based on clustering and the approach presented by Pineda et al. [63]. The procedure is depicted in Figure 5 and also implemented in the open-source tool, tsam. In contrast to the clustering of typical periods, the normalized time series are rearranged in such a way that each time step represents a single sample point, which is, in the case of multiple time series, multi-dimensional as well. Furthermore, the time steps can only be clustered with adjacent time steps in order to form a larger one, i.e., the clustering procedure must respect the side constraint of adjacency. This constraint can be easily respected using hierarchical clustering because of its sequential algorithm, which only merges those (adjacent) clusters that are closest to each other if a so-called connectivity matrix is defined.

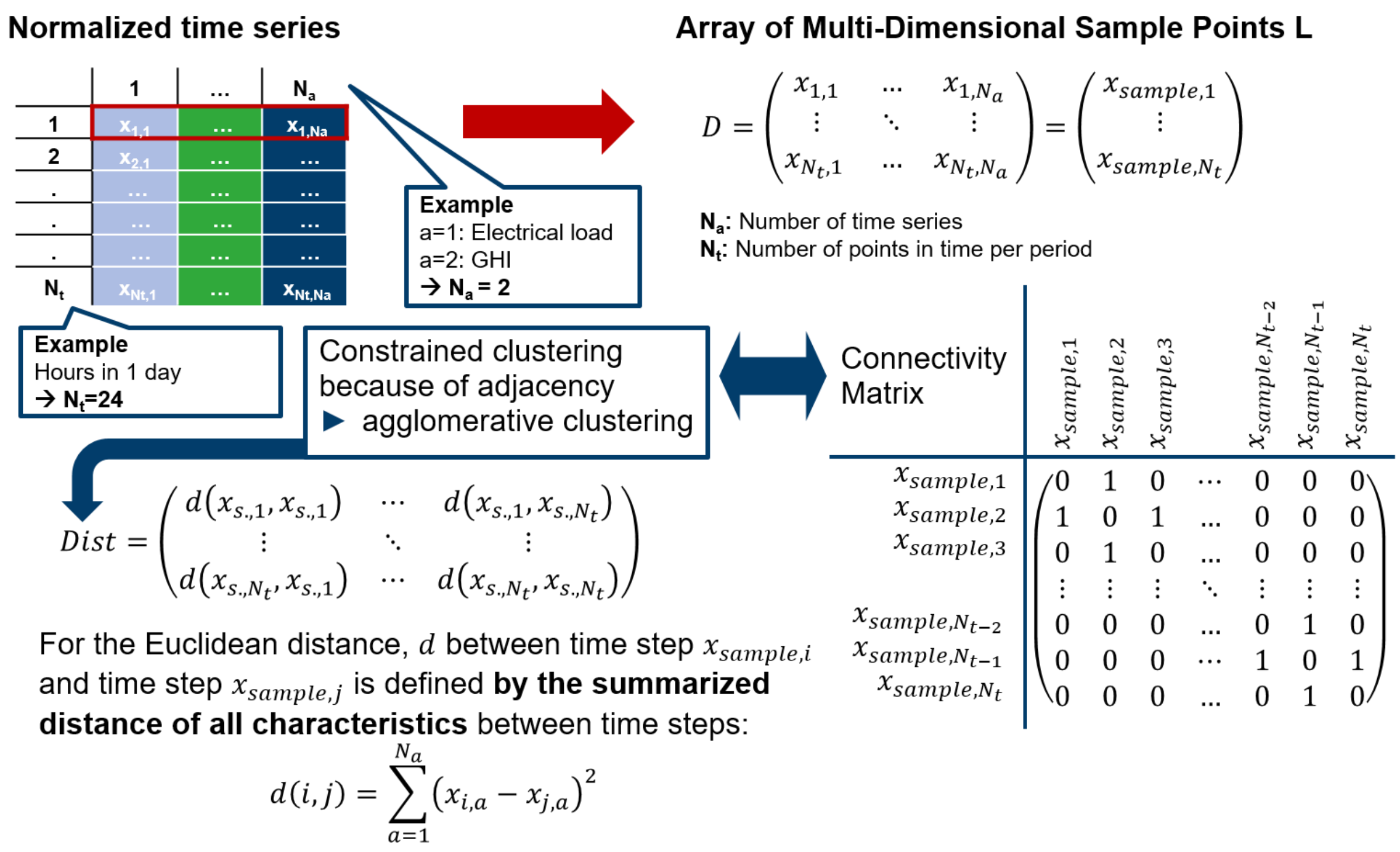

*Figure 5. The procedure of segmentation using constrained hierarchical clustering.*

## 3.2. The optimal tradeoff between segments and typical days

As was shown in the beginning and could be observed in Figure 1, the number of typical periods and segments can be freely combined. Moreover, the number of constraints and variables of an energy system model scales monotonously with the total number of considered time steps, which is defined by the number of periods times the number of segments per period. Therefore, multiple two-tuples of typical period and segment numbers exist with a comparable number of total time steps, and therefore a comparable, expected optimization runtime. Yet, it is obvious that not every one of these configurations is likely to lead to a good approximation of the original time series or a good optimization result. For example, considering 64 typical days with one segment per day instead of eight typical days with eight segments per day will significantly underestimate the intra-day variability, which is crucial for an appropriate representation of technologies, such as photovoltaic panels or daily storage systems.

In order to estimate the quality of aggregation, different accuracy indicators are used in the literature to quantify the difference between the original time series and the modified one. Among these, the Root Mean Squared Error (RMSE) is the most popular, which was either used to quantify the difference between the original and aggregated time series [31, 38], an example of which is depicted in the left half of Figure 6, or between the respective duration curves [15, 31, 67, 68, 72, 81, 85-88], as depicted in the right half of Figure 6.

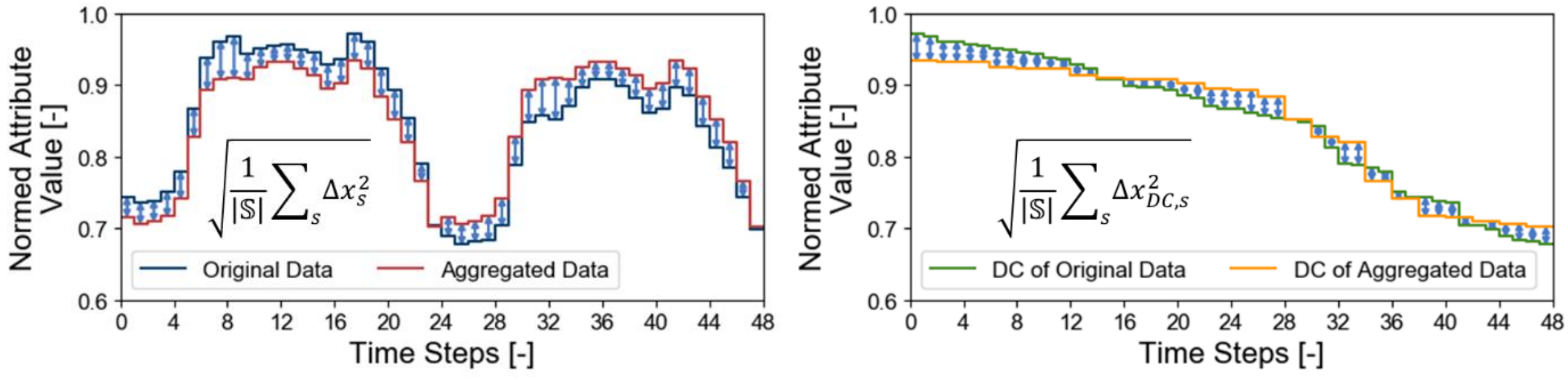

*Figure 6. Illustration of squared summands used for calculating the RMSE (left) and between the respective duration curves (right) as defined in [59].*

In this study, the RMSE between the normed time series and the corresponding aggregated one was used because it takes into account the difference in each time step, rather than the differences in the value distribution of the time series only, which is agnostic to the value chronology in the time series. For multiple time series that are simultaneously aggregated, the RMSE for all of the time series is defined as follows:

$$\mathrm{RMSE_{tot}} = \sqrt{\frac{1}{|\mathbb{A}| \times |\mathbb{T}|} \sum_{|\mathbb{A}|} \sum_{|\mathbb{T}|} \left(x_{a,t}^{aggregated} - x_{a,t}^{original}\right)^2}$$

Here, $|\mathbb{A}|$ is the number of attributes and $|\mathbb{T}|$ the number of time steps.

Indicators such as the RMSE are based on the assumption that a small deviation between the original and the aggregated time series leads to a small deviation between the energy system model's optimization result based on the original time series and that based on the aggregated time series. Moreover, the arbitrary combination of typical periods and segments offer the option of searching for a good tradeoff between them in order to reduce the RMSE as much as possible for a given number of time steps. Therefore, we propose an algorithm based on the method of steepest descent, as suggested by Haskell et al. [89]. In this context, the steepest descent is defined as the ratio of an RMSE's decrease over an increase in the total number of time steps. As the latter can be either increased by increasing the number of typical periods or the number of typical segments, the algorithm checks both options and chooses the more efficient one, an example of which is illustrated for three consecutive steps in Figure 7. The decision rule can be described by the following equation:

$$\hat{t}_{n+1} = \hat{t}_n + \delta \cdot \min\left(\frac{\mathrm{RMSE_{tot}}(\hat{p}_{n+1}, \hat{s}_n) - \mathrm{RMSE_{tot}}(\hat{p}_n, \hat{s}_n)}{\hat{s}_n \cdot (\hat{p}_{n+1} - \hat{p}_n)}, \frac{\mathrm{RMSE_{tot}}(\hat{p}_n, \hat{s}_{n+1}) - \mathrm{RMSE_{tot}}(\hat{p}_n, \hat{s}_n)}{\hat{p}_n \cdot (\hat{s}_{n+1} - \hat{s}_n)}\right)$$

where $\hat{t}_n$ describes the number of total time steps of the aggregated time series at iteration $n$, $\delta$ the finite resolution increase and $\hat{p}_n$ and $\hat{s}_n$the number of typical periods and number of segments at that iteration step, respectively.

The complete algorithm is described in Figure 8. Starting from a single typical period with one segment, i.e., an aggregated time series consisting of a single time step, the algorithm checks in each step what direction of increase in temporal resolution is the most beneficial with respect to a decrease in the RMSE. Once the number of typical periods equals the resolution of the original time series, only increases in segments are considered, or vice versa. Moreover, it is possible to define a maximum number of total time steps for the aggregated time series, as it defines the remaining size of the aggregated energy system model. Accordingly, the algorithm may stop once this total number of time steps has been surpassed. In our calculations, the step width between the number of typical days and segments followed a ratio of approximately $\sqrt{2}$ in each direction, as this proved to be suitable for the degressive RMSE decrease for higher temporal resolutions.

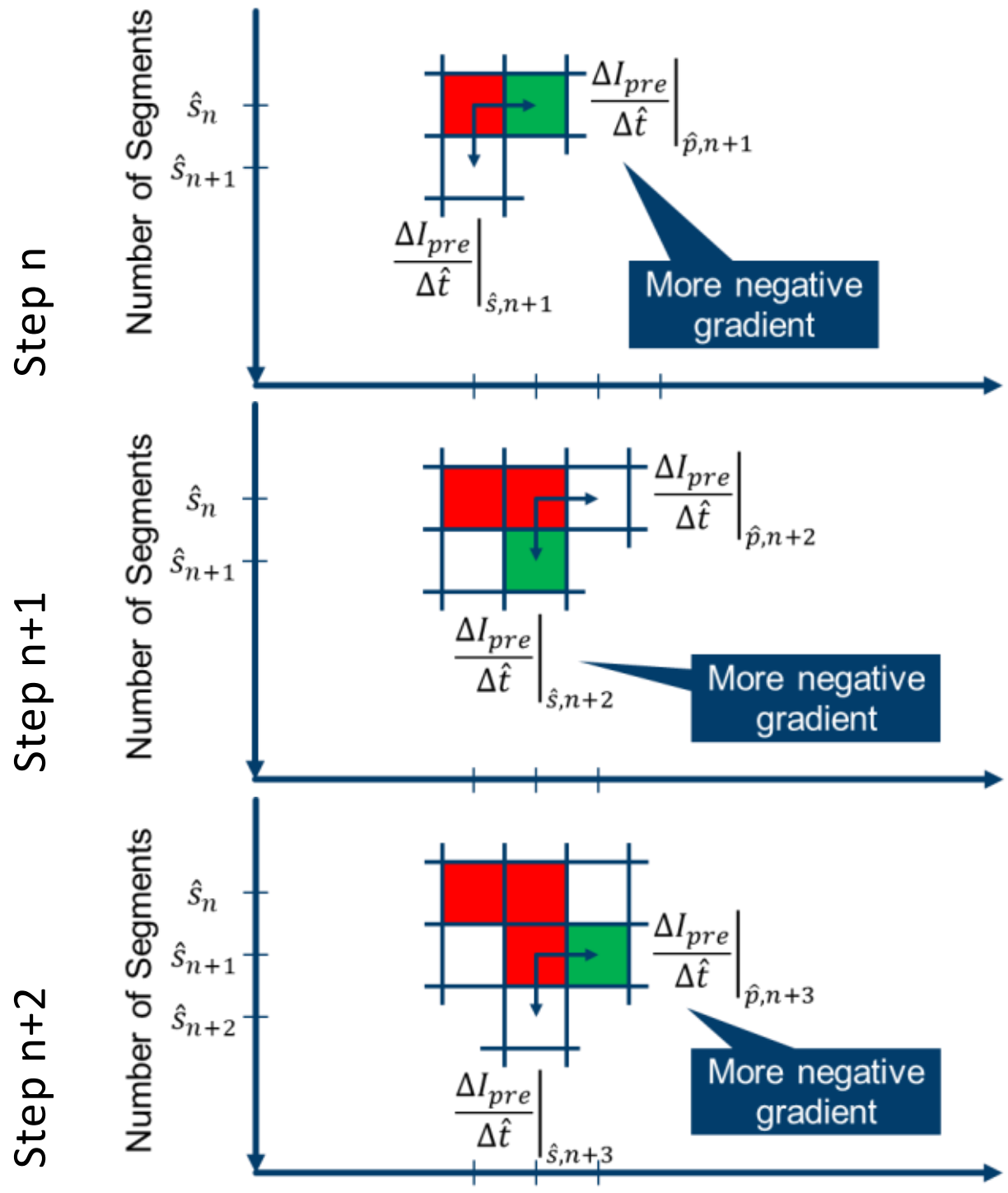

*Figure 7. Three subsequent steps of the proposed algorithm to find the optimal number of segment and typical days.*

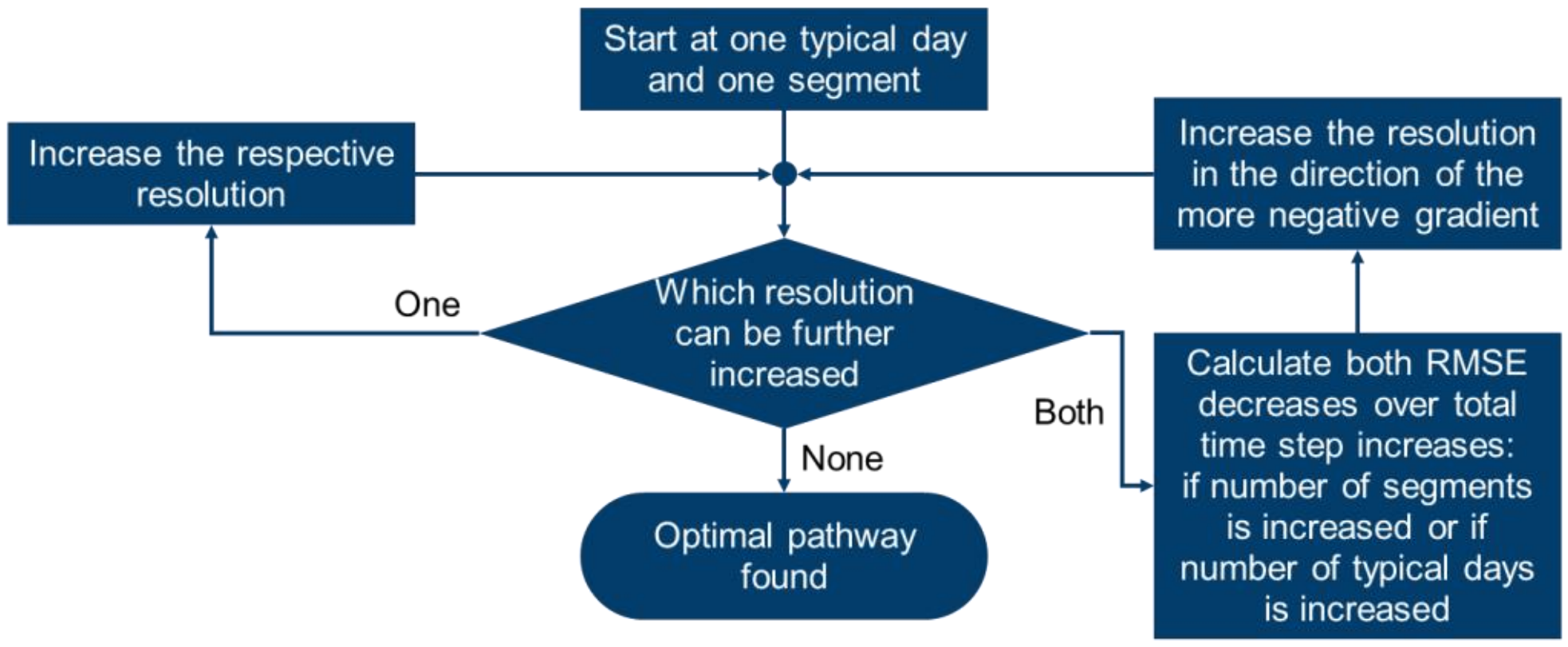

*Figure 8. The optimal pathway algorithm from minimum to maximum resolution.*

As is displayed in Figure 9, the algorithm is capable of differentiating between different types of time series. In the case of a time series for global a horizontal time series, the algorithm prioritizes an increase in the number of segments up to a level of twelve segments per typical days before choosing to increase the number of typical days. In contrast, the wind time series can be aggregated more efficiently by first increasing the number of typical days before the inner-daily resolution is increased. This emphasizes the fact that the algorithm takes the specific features of the time series into account. While solar irradiance profiles exhibit a clear daily pattern, are close to zero during the nighttime hours, and have steep gradients during the morning and evening hours, wind time series feature an aperiodic pattern with fairly smooth hourly transitions. Therefore, solar irradiance profiles require a high intra-daily temporal resolution to appropriately cover about 12–16 hours of the day. During the nighttime hours, a significantly smaller temporal resolution is sufficient. In contrast, wind speed time series require many typical days due to their aperiodic behavior, leading to possible wind peaks at any point in the daytime, whereas a relatively small temporal resolution within each day is sufficient due to the relatively slow transition in the wind speed. In an energy system optimization, multiple time series are considered together. Therefore, the observed examples constitute extreme cases. In general, the algorithm adapts to the specific types of considered time series due to the definition of the $RMSE_{tot}$, which enables the simultaneous evaluation of all time series.

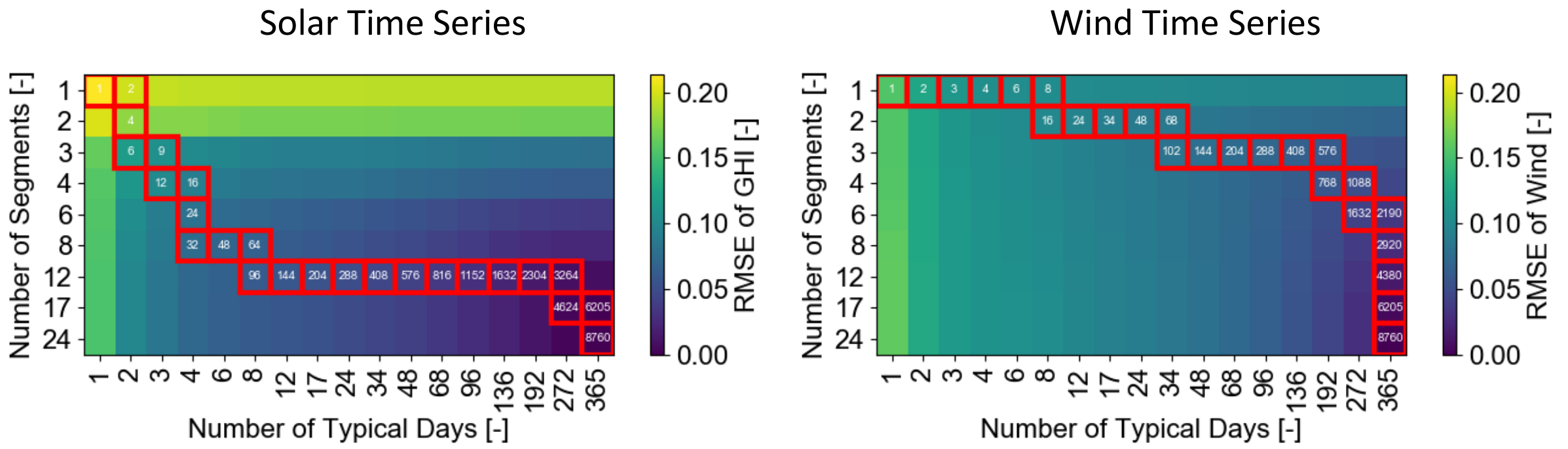

*Figure 9. Proposed algorithm for the RMSE of a typical global horizontal irradiance and a wind speed profile.*

## 3.3. Preserving the original time series' value distribution

The second algorithm proposed in this work in order to improve the quality of temporal aggregation for energy system models addresses the underestimation of extreme values as observed for k-means and k-medoids. Although numerous studies have manually included extreme values from the original time series into the aggregated ones, this procedure requires profound knowledge regarding the potentially critical system situations and suffers from the combinatorial problem that the number of extreme values grows in conjunction with the number of time series [8]. Therefore, the algorithm proposed below strives for the precise reproduction of the original time series' value distributions (value curves).

The algorithm consists of the following steps:

1. For a specific attribute $a$, sort all values of those periods that are assigned to a specific cluster $k$ and yield the cluster's duration curve $\widetilde{\mathbb{X}}_{a,k}$. This is illustrated by the green duration curve shown in Figure 10.
2. Average every $|\mathbb{C}_k|$ values of the cluster's duration curve $\widetilde{\mathbb{X}}_{a,k}$ in order to obtain the duration curve of the cluster's representative $\widetilde{\mathbb{Y}}_{a,k}$, which is represented by the brown line in Figure 10.
3. At the same time, calculate the mean (i.e., centroid) profile $\mathbb{M}_{a,k}$ for the specific attribute $a$ using all periods assigned to cluster $k$ for each time step represented by the dark blue line in Figure 10.
4. Determine the mean profile's duration curve $\widetilde{\mathbb{M}}_{a,k}$ by sorting its values and extract the order of time steps $t$ in which they appear in the sorted curve $\widetilde{\mathbb{M}}_{a,k}$. The sorted mean profile is represented by the yellow line in Figure 10.
5. Assign the index order from $\widetilde{\mathbb{M}}_{a,k}$ to the sorted values of $\widetilde{\mathbb{Y}}_{a,k}$ and sort them such that the attached indices of $\widetilde{\mathbb{Y}}_{a,k}$ are in ascending order. The result is the representative profile $\mathbb{Y}_{a,k}$ and is depicted as a red line in Figure 10.

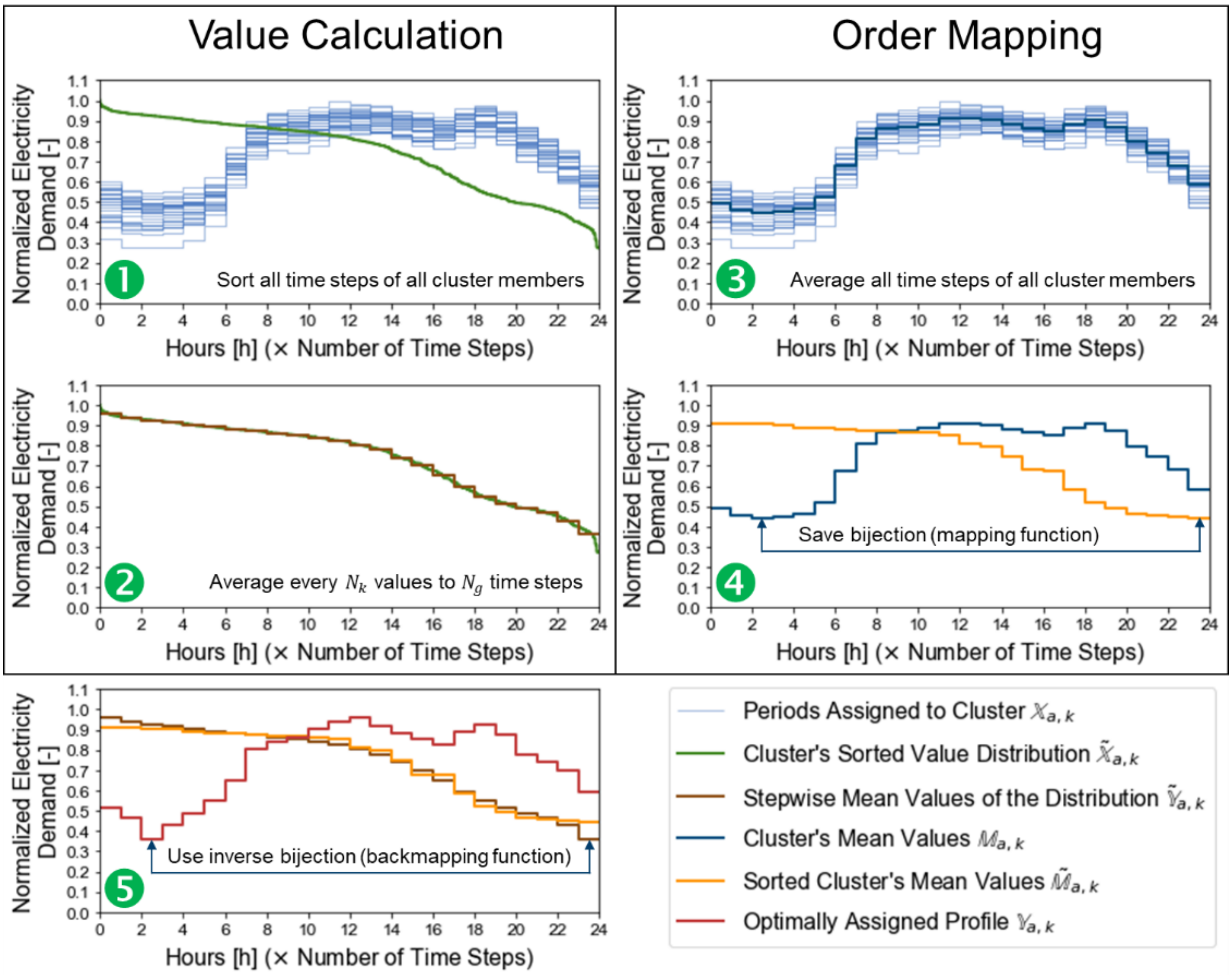

*Figure 10. Graphical interpretation of the algorithm.*

In contrast to a cluster representation by centroids or medoids, this approach is based on a direct approximation of the original time series' value distribution. The impact of this approach on a typical electricity time series with hourly resolution can be observed in Figure 11. On the left-hand side, the duration curve of the original and aggregated time series, consisting of eight typical days at hourly resolutions, are depicted. It can be observed that the original time series' duration curve is almost perfectly replicated. The right-hand side of Figure 11 depicts the original and aggregated electricity profiles during nine days in February. As the proposed method preserves the original time series' value distributions and therefore approximates the extreme values as well, the aggregated profile covers the cluster's extreme periods within its respective representative, which leads to a more peaky profile.

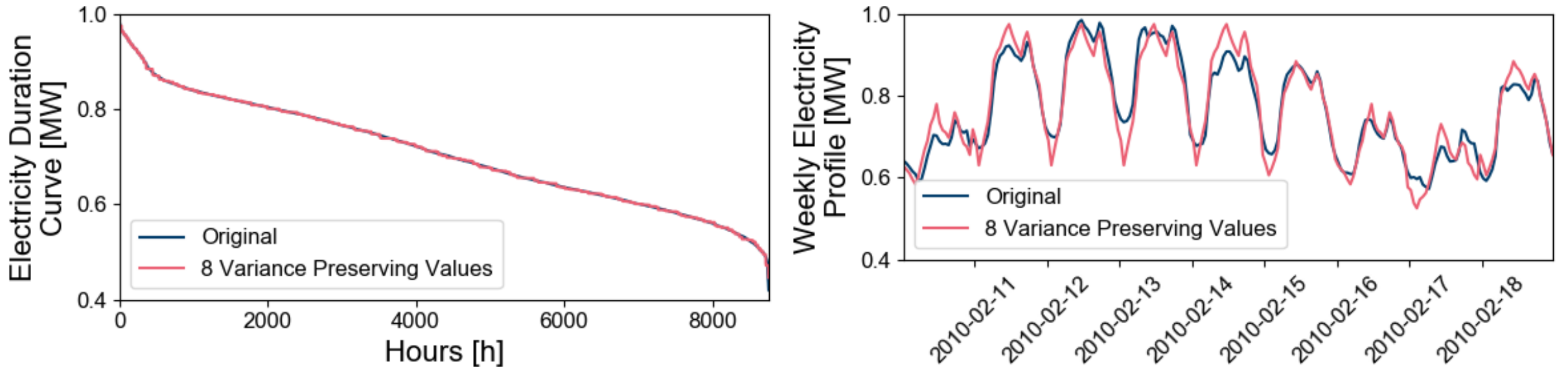

*Figure 11. Yearly duration curve of a typical electricity profile, its profile during a week in February, and the corresponding aggregated time series as predicted by eight typical days using hierarchical clustering and the approximate distribution-preserving representation method*

## 3.4. Case studies

In order to validate the effectiveness of our proposed methods, we introduce two different capacity expansion models that differ significantly with respect to their number of regions, their mathematical model structure, and their considered technologies. Both models have been implemented in the FINE[3] optimization framework [44].

### 3.4.1. The self-sufficient building model

The first case study used for validation was introduced by Kotzur et al. [90] and is a model for a self-sufficient building that relies exclusively on photovoltaic electricity feed-in, which is also used to cover household heat demand using a heat pump and an electric boiler. In order to account for daily and seasonal variations in the electricity supply, the system comprises a battery, a heat and two hydrogen storages, as depicted in Figure 12. The hydrogen subsystem is linked to the heat and electricity system by a reversible solid-oxide fuel cell (rSOC). Furthermore, the system has been used as an application case in earlier studies [59, 91] and its techno-economic data is provided in Appendix A.1, and was made publicly available in [92].

The system is a single-node model considering two time series for heat and electricity demand, as well as two rooftop photovoltaic and one open field photovoltaic time series generated using the open-source software tsib [93, 94]. For this, the weather year of 2013 in the location of Berlin was taken from the COSMO dataset [95-97]. Additionally, the energy system model considers discrete investment decisions for certain components using binary variables, turning the model into an MILP.

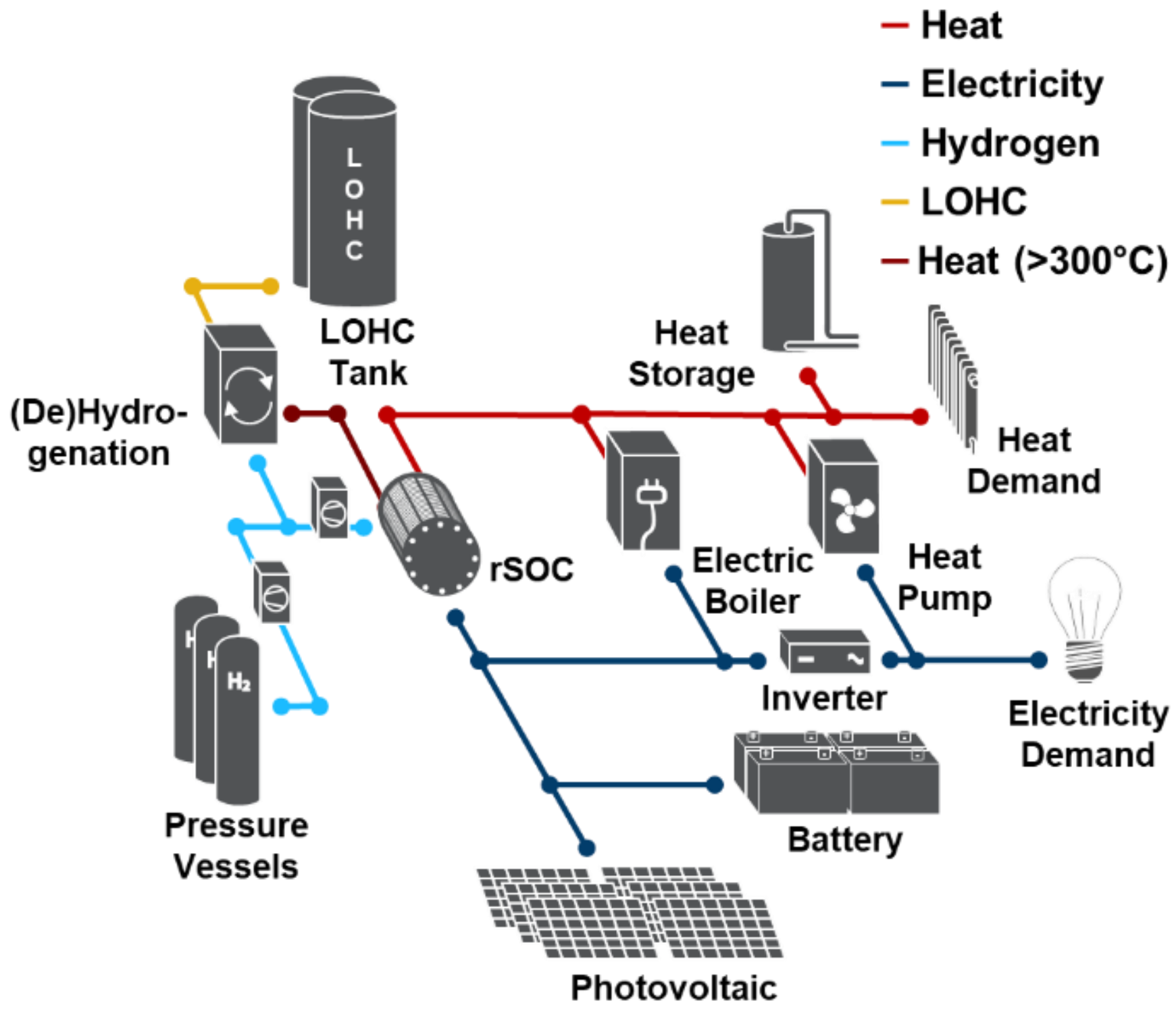

*Figure 12. Simplified scheme of the self-sufficient building model.*

[3] FINE package available at: https://github.com/FZJ-IEK3-VSA/FINE (date accessed: 11/23/2021)

### 3.4.2. The European model

In contrast to the self-sufficient building model, the European model developed by Çağlayan et al. [98, 99] is a carbon-neutral capacity expansion model with the target year 2050, and considers 96 European regions. Moreover, it is not a greenfield study and accounts for preexisting infrastructure.

As is shown in Figure 13, the system primarily focuses on a hydrogen and electricity network based on renewable energy sources, such as PV, wind energy, biomass, and hydroelectricity. The hydrogen and electricity networks are connected by electrolyzers, gas engines, and different types of fuel cells and gas turbines. Furthermore, hydrogen and battery storage is considered in order to account for fluctuating supply and demand. The hydrogen demand is assumed to originate exclusively from fuel cell-electric vehicles with a market penetration of 75% [100], whereas the electricity demand was obtained from the E-Highway study [101] for a 100% renewable scenario. This demand is considered to account for electrified heating in the residential and industrial sectors, as well as the operation of battery-electric and plug-in hybrid vehicles.

For the fully resolved case, the European model is on the edge of computational feasibility, with a calculation time of approximately four days, and therefore does not consider binary decision variables, i.e., it is a large LP. Furthermore, the multitude of regions and technologies leads to more than 900 time series, which are then simultaneously clustered. These constitute fundamental differences from the self-sufficient building model with respect to both the application case and mathematical structure of the model.

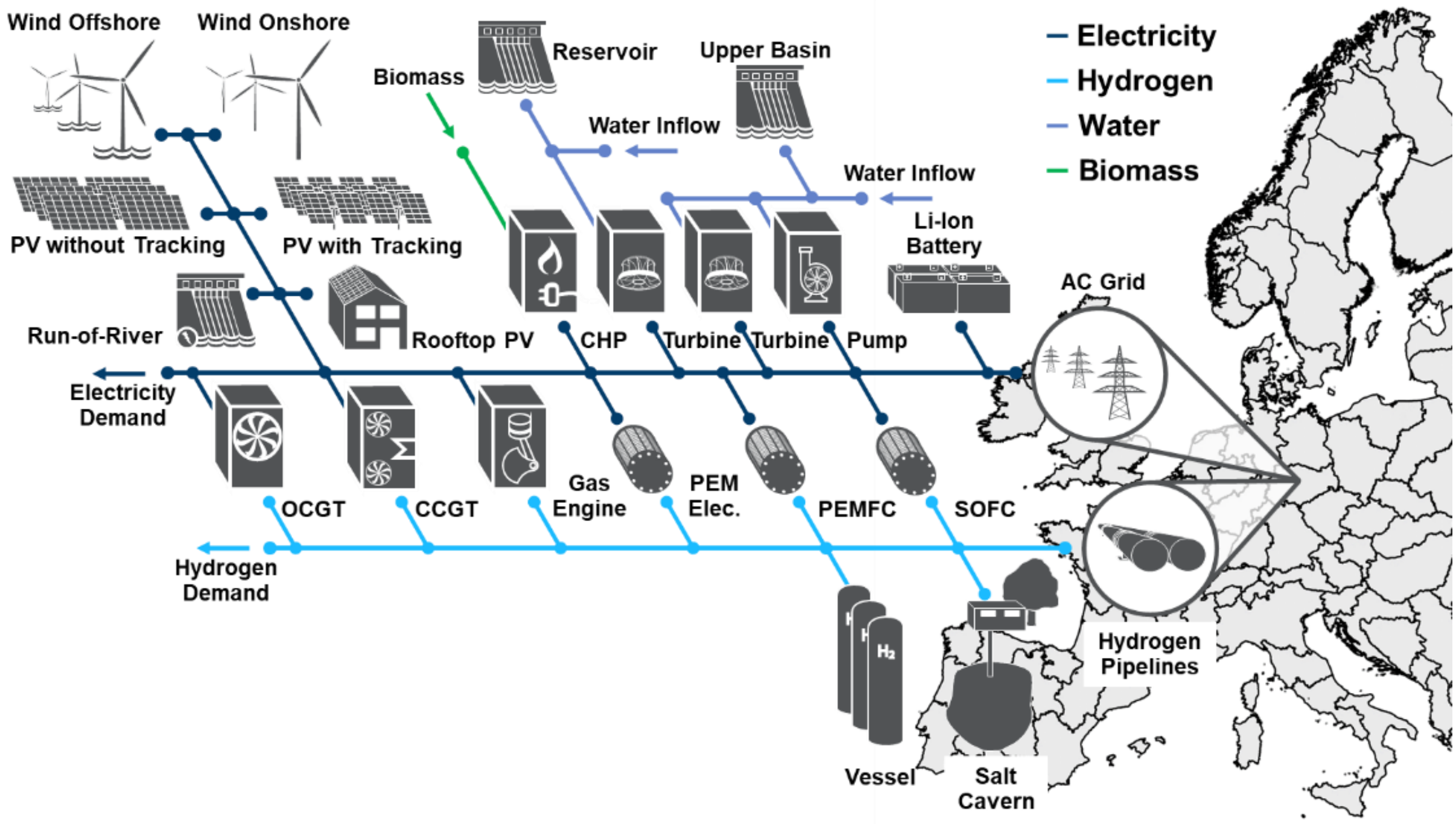

*Figure 13. Layout of the European model.*

## 3.5. Computational resources and comparison framework

As previously noted, the computational runtime not only depends on the size of the model but also on the hardware used. For that reason, the computational resources used to conduct the analyses in Section 4 are also listed in Table 3. All calculations were performed on a computer cluster in parallel, but with an individual core assigned to each calculation. Gurobi 9.0.0 was used to carry out the optimization, with the thread number set to one in order to prevent unpredictable variations in the calculation times arising due to data traffic between the threads.

*Table 3. Computational resources used for Section 4.*

| | |
|---|---|
| CPU Model | Intel(R) Xeon(R) Gold 6144 CPU |
| Number of Cores per Computing Node | 16 |
| Threads per Core | 2 |
| CPU Max Frequency [MHz] | 3.5 |
| Shared Memory [GB] | 1024 |

## 4. Results

In this section, the models introduced in Section 3.4 are optimized for the large set of different aggregation configurations depicted in Figure 14. Each field describes an aggregation setup, with the white number representing the total number of remaining time steps for that configuration, .e.g., in the case of eight typical days and eight segments, the time series were aggregated from 8760 to 64 time steps. The ratio between two neighboring typical period or segment numbers was selected to be approximately $\sqrt{2}$, as this accounts for a degressive trend of accuracy over an increase in the total number of time steps. Moreover, only typical periods with a length of one day were considered due to the daily charging pattern of short-term storage in the respective models. A more detailed discussion on the optimal period length for temporal aggregation is provided in an earlier work by the authors [59].

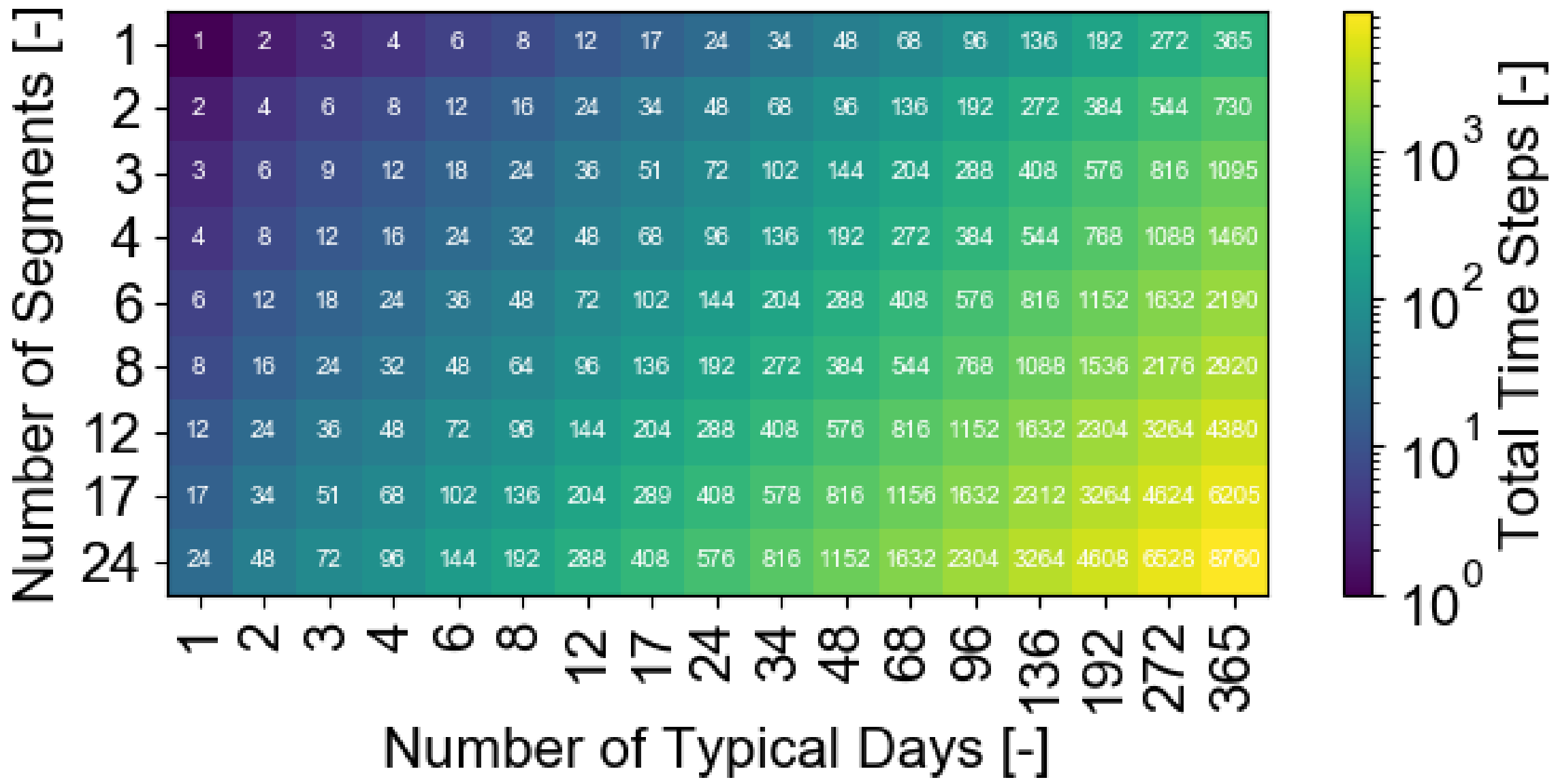

*Figure 14. Number of total time steps depending on typical day and segment configurations that were considered in the sensitivity analysis for four different representation methods*

Furthermore, each typical day and segment configuration was represented by centroids, medoids, and the proposed distribution-preserving algorithm, leading to 9 x 17 x 3 = 459 different optimization runs per model. Throughout the entire chapter, Ward's hierarchical clustering algorithm was used due to its good runtime scaling behavior with respect to both the sample number and sample dimensionality. For the considered multi-regional model, it was the only computationally-tractable clustering algorithm and, for the single-regional model, it was found that the results were less sensitive to the clustering technique than the representation method. In addition, the reproducibility of the results was given due to its deterministic algorithm.

After all optimizations were run, the most commonly used aggregation configurations from the literature, i.e., typical days at hourly resolution, represented by either centroids or medoids, were benchmarked against the optimization results obtained from the aggregation configurations, along with the optimal ratio of typical days and segments using the distribution-preserving algorithm as the representation method.

## 4.1. Results for the self-sufficient building

First, the large-scale parameter variation was applied to the self-sufficient building model. The right-hand side of Figure 15 depicts the individual model runtimes and relative cost deviations of the respective aggregated models from the fully resolved reference case for all aggregation configurations. Here, the chosen representation methods are highlighted in different colors.

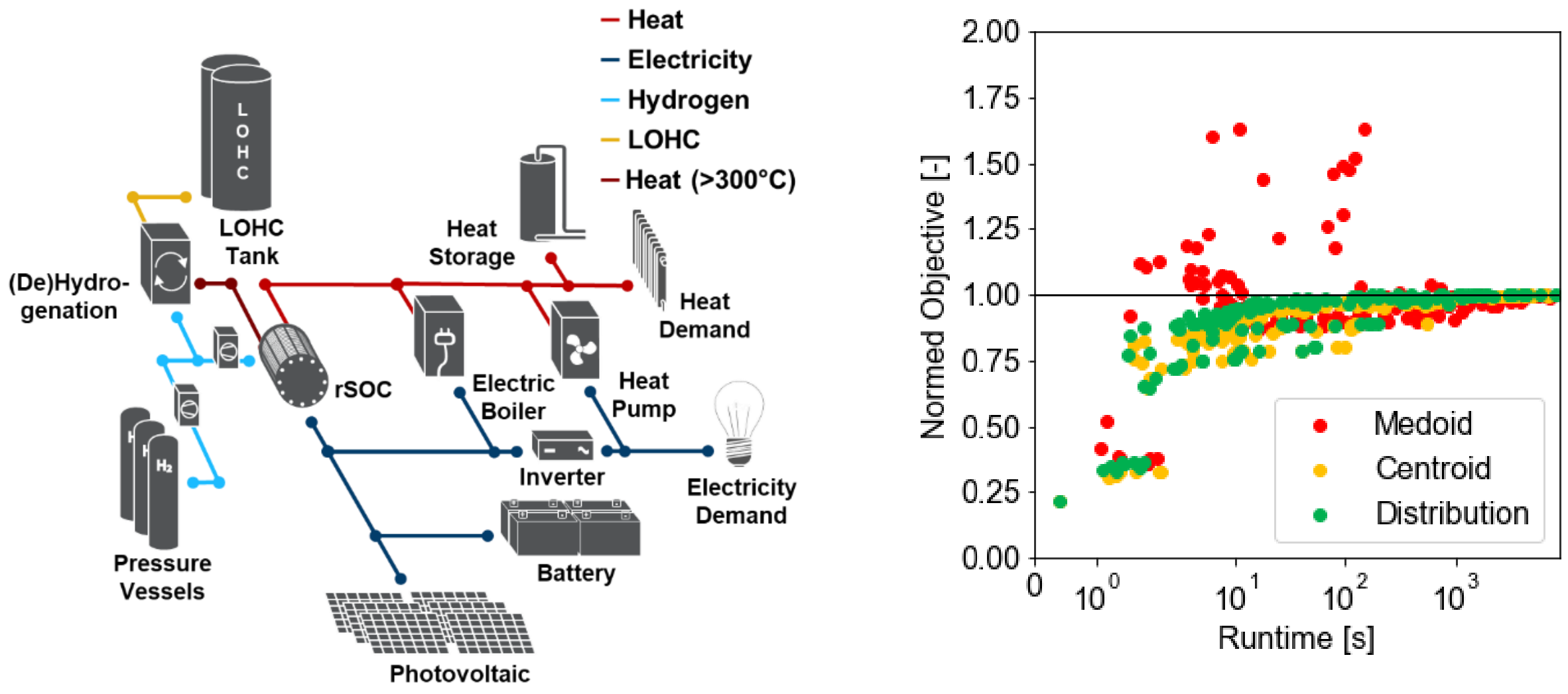

*Figure 15. The self-sufficient building model and corresponding normed optimal objectives over the computation times for all considered configurations.*

As Figure 15 indicates, there is a clear tradeoff between the runtime of the models and their respective accuracies, which is represented by the aggregation-induced deviation from the optimal objective function values of the reference case. The more time steps that are considered, the more the system costs approach those of the reference case. Moreover, it can be seen that a representation by centroids and the distribution-preserving algorithm lead to a consistent underestimation of the total system costs, but converge somewhat monotonously on the system costs of the reference case for higher resolved temporal aggregation configurations with longer runtimes. In contrast, this trend cannot be observed for a representation by medoids, which may lead to either an over- or underestimation of the total system costs.

The left-hand column of Figure 16 depicts the pathways of the aggregation configurations as proposed by the algorithm introduced in Section 3.2. It is worth noting that the clusters as defined by Ward's hierarchical algorithm for segments and typical days are identical for all three representation methods. However, the different representation methods lead to slightly differing aggregated time series, and therefore different RMSE values and proposed ratios between the typical days and segments. For all representation methods, it can be clearly seen that the algorithm initially favors a relatively high number of segments for relatively strongly aggregated time series. For twelve typical days, the proposed number of segments is between six and eight. However, for larger numbers of total time steps, the number of necessary segments per typical day is at most twelve, except for the configurations with 365 typical days, which is equal to the original number of days. Both observations can be explained by the heavy reliance of the self-sufficient building model on PV, for which this effect was already explained in Section 3.2. Moreover, the remaining time series pertain to heat and electricity, which also have a strong daily pattern and fairly smooth profile intervals between midnight and 6:00 AM. Accordingly, the segments can be comparably large during the night, which leads to a sufficient model aggregation with no more than twelve time steps per day.

The graphs shown in the right half of Figure 16 map the aggregation-induced deviations of the optimal objective onto the respective aggregation configurations. Here, the convergent behavior for a representation by centroids or the distribution-preserving algorithm can again be observed, indicating that a higher number of time steps indeed leads to smaller deviations. Likewise, it can also be seen that an increase in either the number of typical days or segments per typical day always drives the need to increase the complement resolution, e.g., in the case of four typical days represented by centroids, it does not make a big difference whether four segments per day are chosen or 24. If the total “runtime budget” would not allow for more than 96 total time steps, a better choice would therefore be to choose twelve typical days with eight segments each instead of four typical days with 24 each. For a representation by medoids, however, an unpredictable behavioral pattern can be observed that impedes an intuitive choice for a sufficient number of typical days or segments. Assuming that a system is not solvable in its fully resolved configuration, a sensitivity analysis in the neighborhood of a specific temporal aggregation configuration would therefore not lead to easily interpretable results, and is therefore not reliable.

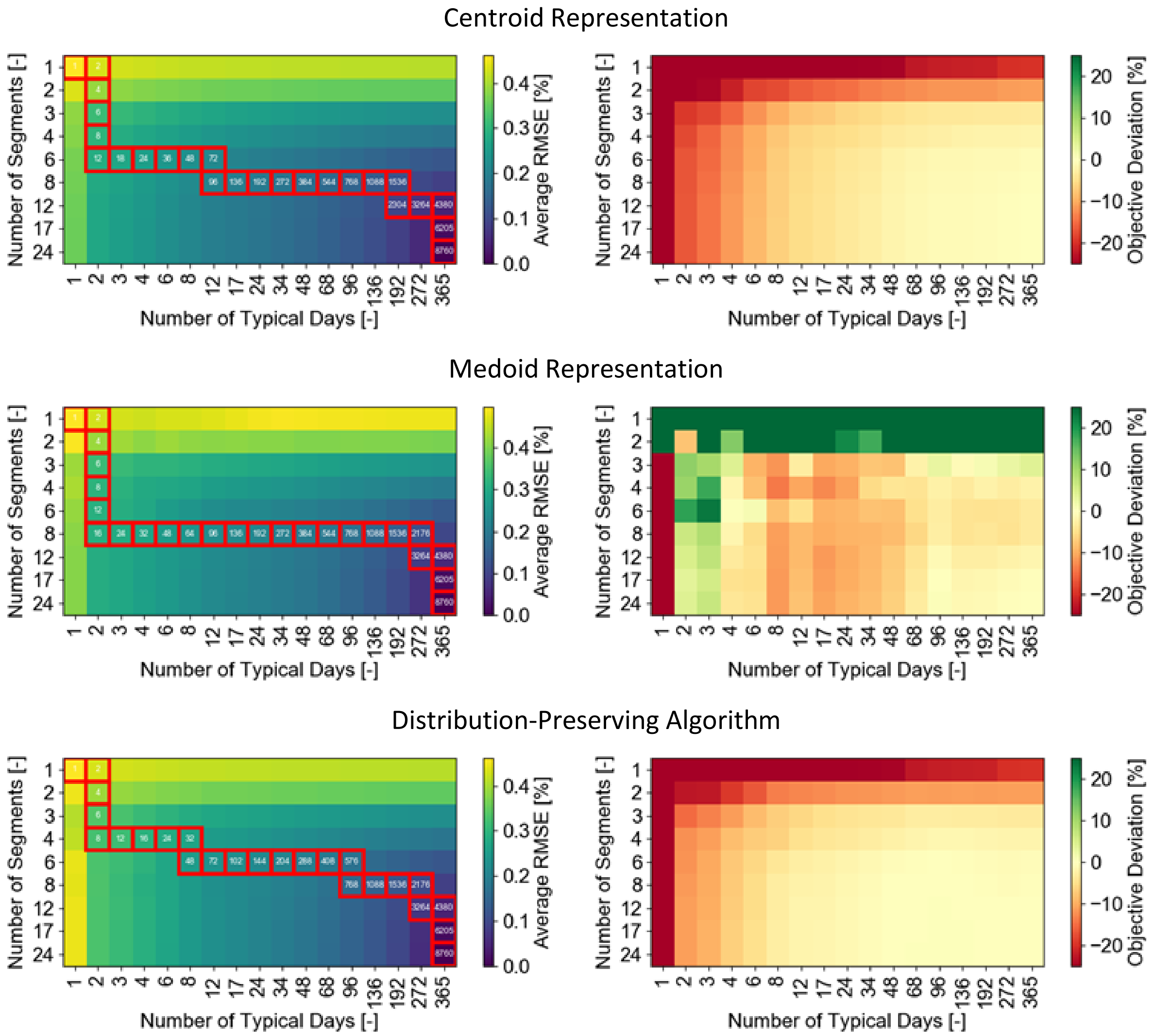

*Figure 16. The pathway found by the proposed algorithm for an optimal ratio between the number of typical days and segments depending on the representation method (left) and deviation from the optimal objective of the fully resolved case depending on the configuration of the number of typical days and segments, as well as the representation method (right).*

In order to evaluate the effectiveness of choosing a balanced number of typical days and segments for each typical day, Figure 17 separately depicts the results shown in the right half of Figure 15 for each representation method. Moreover, the aggregation configurations that are most frequently used in the literature, i.e., typical days with an hourly resolution, are highlighted as violet dots, whereas the configurations as obtained by an improved ratio between the number of typical days and the number of segments are highlighted as green dots.

As medoids do not lead to consistent results that monotonously converge upon the reference case, it is self-explanatory that choosing a balanced number of typical days and segments within each typical day does not lead to predictable results either. However, those configurations requiring more than 10 seconds of optimization runtimes are indeed closer to the objective functional value of the reference case than those obtained by typical days with an hourly resolution. In contrast, the algorithm proposed in Section 3.2 not only outperforms the configurations of typical days with an hourly resolution but also leads to Pareto-optimal results in the case of a representation by centroids or when using the distribution-preserving algorithm.

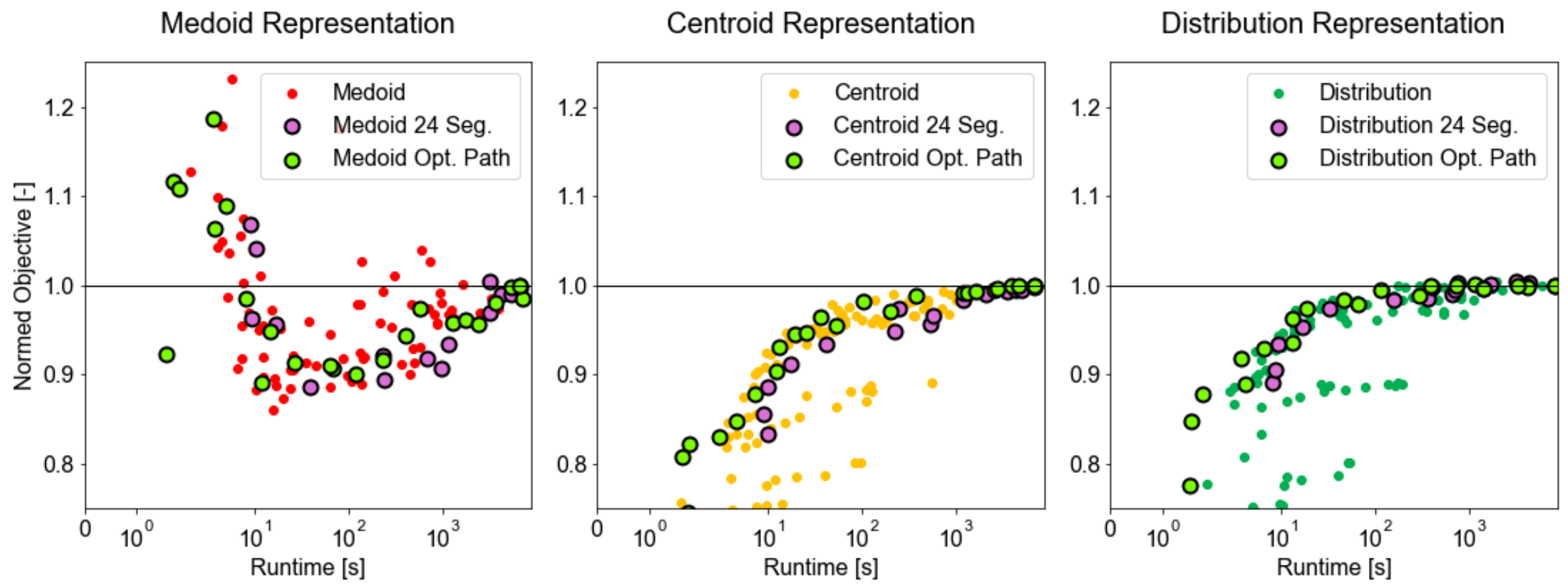

*Figure 17. Convergence behavior of the different representation methods if an optimal ratio between the number of typical days and the number of segments is chosen based on the RMSE compared to the convergence behavior of different typical day numbers at hourly resolutions.*

Lastly, the combination of both of the proposed improvements on current aggregation approaches is compared to all other model runs in the left graph of Figure 18. Here, it can be seen that those aggregation configurations that were obtained using the proposed number of typical days and segments, as well as the distribution-preserving algorithm, leads to the Pareto-optimal model aggregations. The left graph of Figure 18 compares the same aggregation configurations to the status quo of most aggregation approaches, i.e., a representation by centroids or medoids of typical days with hourly resolutions, which are clearly dominated. Compared to a representation by centroids, it can be seen that the proposed aggregation configurations lead to less than half of the deviations from the optimal objective of the reference case at comparable runtimes. Because of the decreasing gain in accuracy for longer runtimes, this means that the proposed algorithm can obtain results of comparable accuracy in about one tenth the runtime required for representation by centroids.

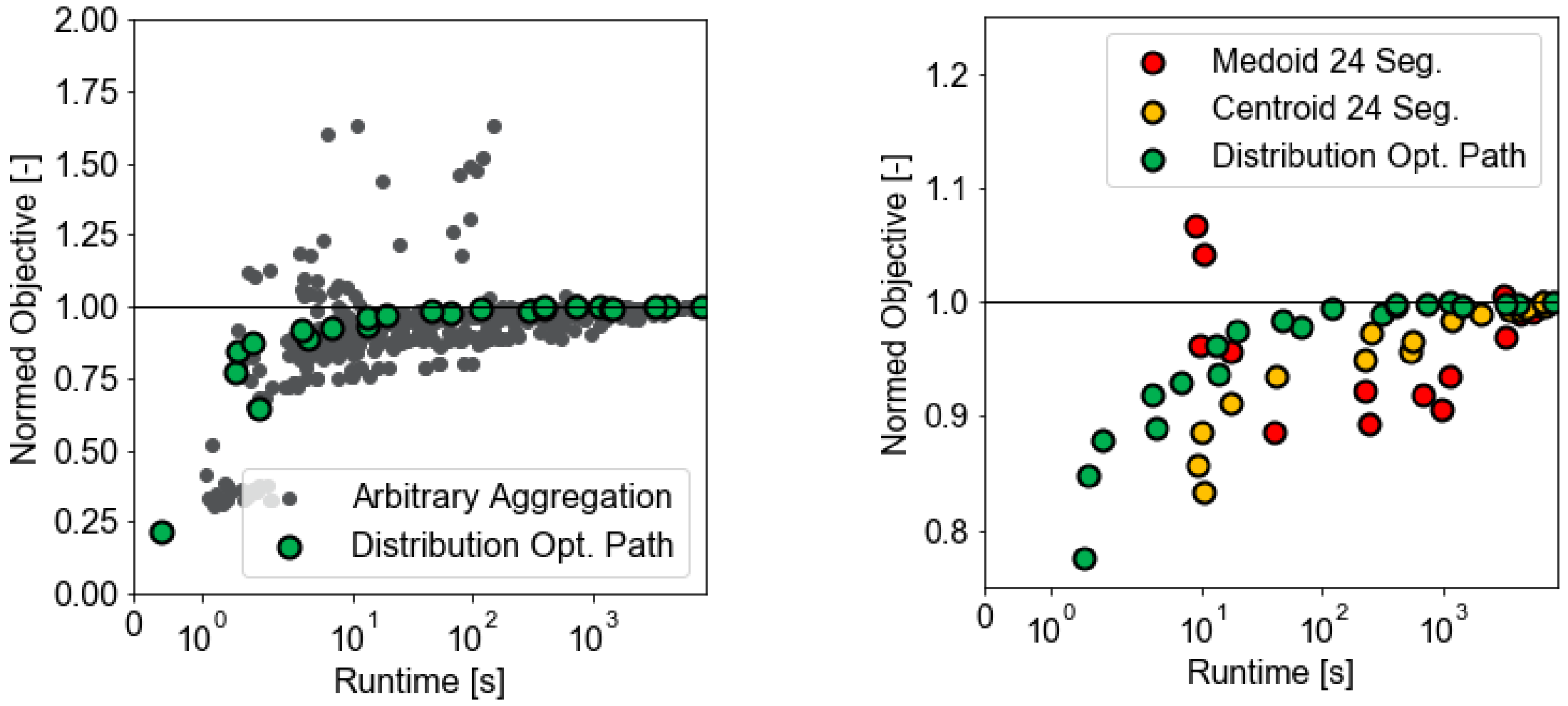

*Figure 18. The identified optimal aggregation method amongst all tested aggregation approaches (left) and a direct comparison to state-of-the-art aggregation methods as used in **FINE** (right).*

## 4.2. Results for the European model

The European system model differs considerably from the self-sufficient building one with respect to the number of time series, the number of regions, and the fact that numerous aperiodic wind time series are considered. Therefore, it can be expected that the results of the different aggregation methods significantly differ. Analogously to the previous section, the right half of Figure 19 depicts the normed objective of all model runs over their respective runtimes separated by the representation method. Again, the representation by centroids leads to a significant underestimation of the system costs, whereas the distribution-preserving algorithm quickly converges on the objective function value of the reference case for increasing numbers of time steps and higher runtimes. In contrast to the self-sufficient building, however, the representation by medoids leads to far less unpredictable behavior. Although medoids still lead to the highest approximations of total system costs, the range of medoid-based optimization runs resembles those of the representation by centroids or the distribution-preserving algorithm. One reason for this is the fact that the representation by medoids does not preserve the original time series' mean values, and can therefore distort the cumulative energy consumption of small island systems, such as the self-sufficient building. For large-scale problems such as the European model, however, the selection of real days from the original dataset accounts for real correlations between time series, and therefore has advantages for large and complex datasets.

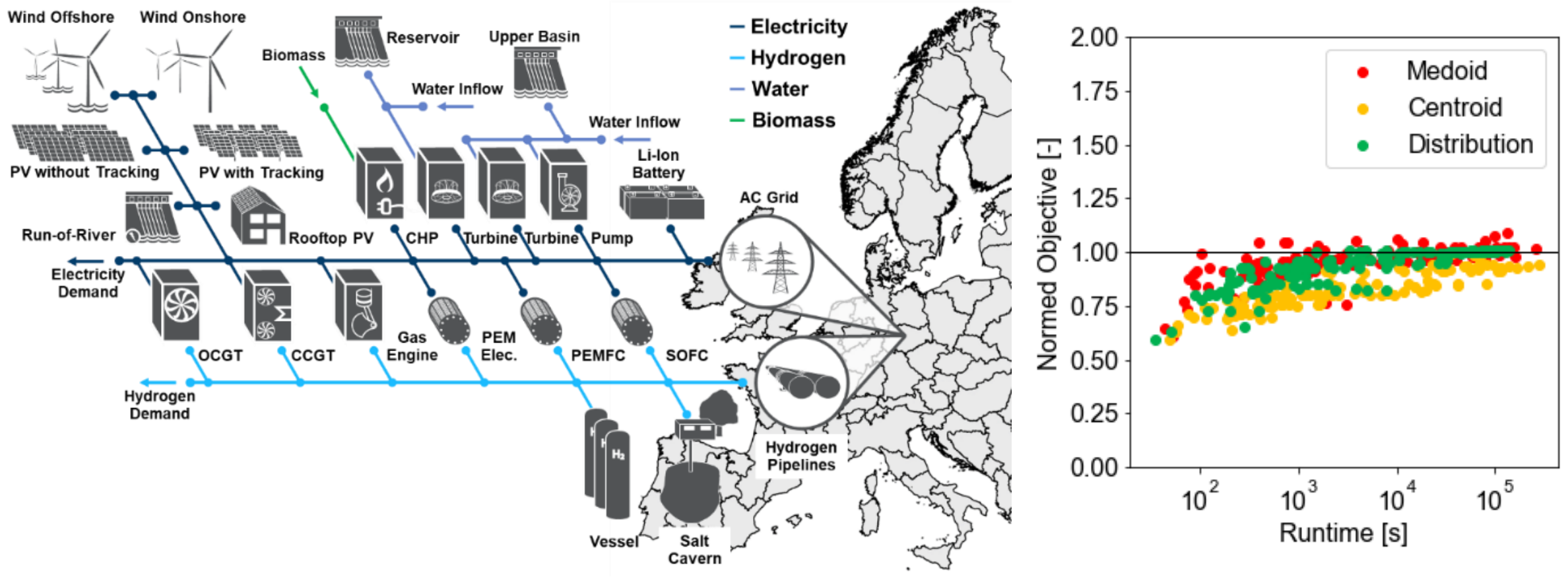

*Figure 19. The European model and corresponding normed optimal objectives over computation times for all considered configurations.*

The left column of Figure 20 depicts the configurations of typical days and segments as proposed by the algorithm presented in Section 3.2. The white fields are optimization runs that were terminated due to numerical instabilities, infeasibility, and non-convergence within the runtime limit set to five days. Here, it can be seen that the algorithm opts to increase the number of typical days first and after that, the number of typical days and segments in a more balanced way than for the self-sufficient building, due to the simultaneous consideration of time series with and without a daily pattern.

The right column of Figure 20 depicts the deviations of the objective function values depending on the aggregation configuration. Again, it can be seen that the total annual costs are relatively indifferent to the number of segments if at least three of them per typical day are chosen. Additionally, medoids do not exhibit a monotonous convergence behavior for higher resolved configurations, and the distribution-preserving algorithm provides a close approximation of the reference case if at least 34 typical days with three segments per day are chosen, which equals a time step reduction of almost 99% compared to the 8760 original time steps, and outperforms all other representation methods.

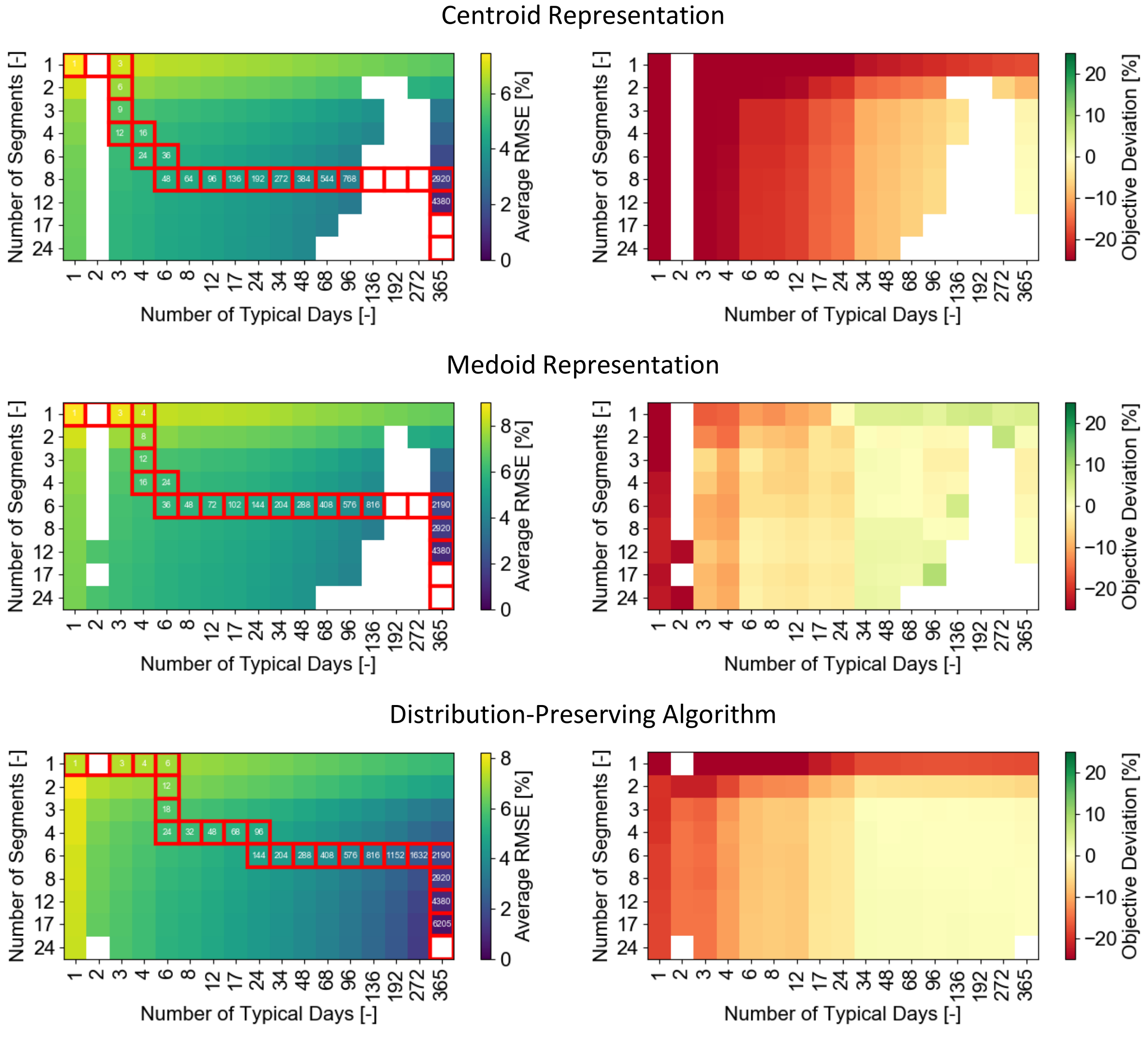

*Figure 20. Pathway found by the proposed algorithm for an optimal ratio between the number of typical days and segments depending on the representation method (left) and deviation from the optimal objective of the fully-resolved case depending on the configuration of the number of typical days and segments, as well as the representation method (right).*

Figure 21 compares the results of the aggregation configurations as proposed by the algorithm in Section 3.2 to those obtained by using typical days with an hourly resolution only, as well as all remaining aggregation configurations. Although the approach does not provide the Pareto-optimal configurations in the case of medoids and centroids, it clearly outperforms the configurations considering typical days at an hourly resolution only. It is worth noting that the green dot lying significantly over the normed objective function value in the case of medoids was a suboptimally-terminated optimization run, more precisely the run with 136 typical days and six segments per typical day, and can therefore be neglected. These findings emphasize that the algorithm proposed in Section 3.2 is agnostic to the financial impact of the individual time series types, as they are an outcome of the optimization itself and therefore do not necessarily choose the optimal configurations from an output data perspective, i.e., with respect to the deviation of the objective function value. However, this should not be the only quality criterion, as a small overall cost deviation does not necessarily indicate meaningful technology sizing, as will be demonstrated in Section 4.3.

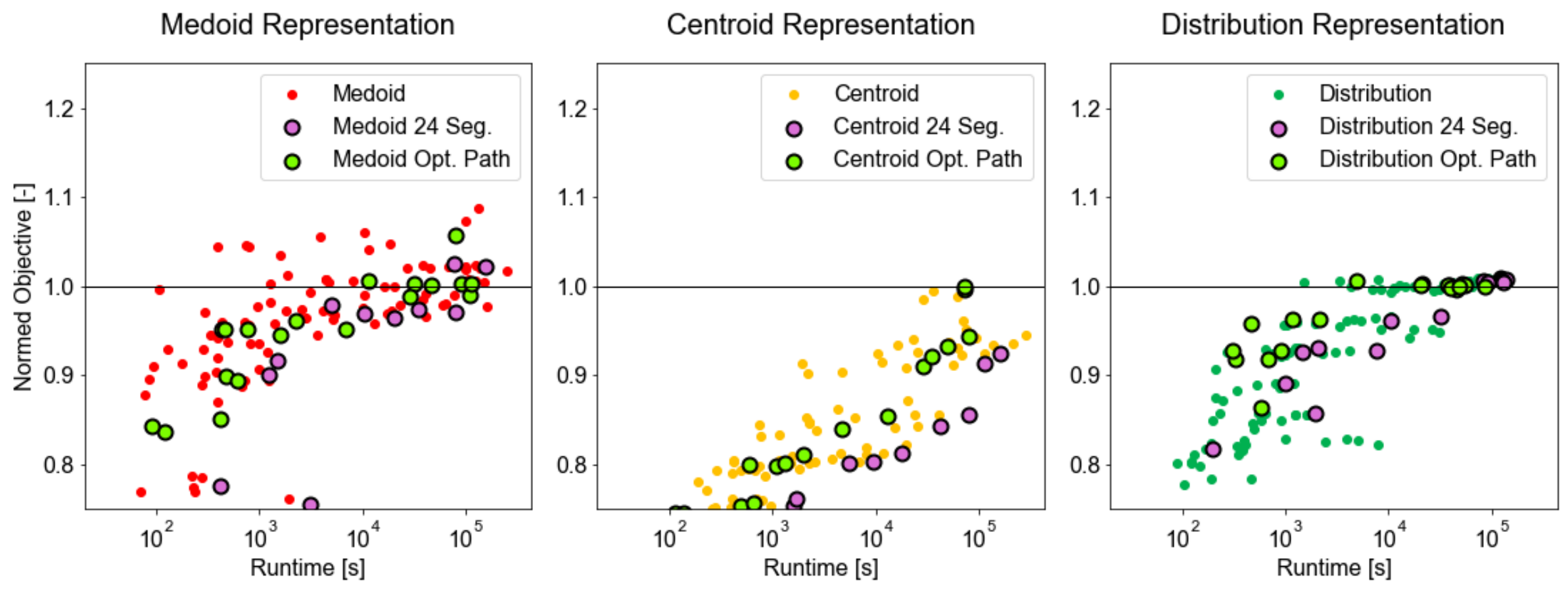

*Figure 21. Convergence behavior of the different representation methods if an optimal ratio between the number of typical days and the number of segments is chosen based on the RMSE compared to the convergence behavior of different typical day numbers at an hourly resolution.*

Finally, Figure 22 visualizes the impact of the proposed improvements compared to all tested aggregation configurations and the status quo of centroid- or medoid-based typical days at hourly resolutions. Here, it can be seen that the proposed methods lead again to the Pareto-optimal solutions for those configurations, with a runtime of more than 1000 seconds. For shorter runtimes and therefore stronger aggregated configurations, the result is not as clear because the distribution-preserving algorithm is based on synthesized time step values, and in the case of too strong an aggregation, the creation of a few typical days with the value distribution of the original time series may distort the correlations between them. In summary, however, the proposed improvements significantly outperform the status quo of the aggregation approaches.

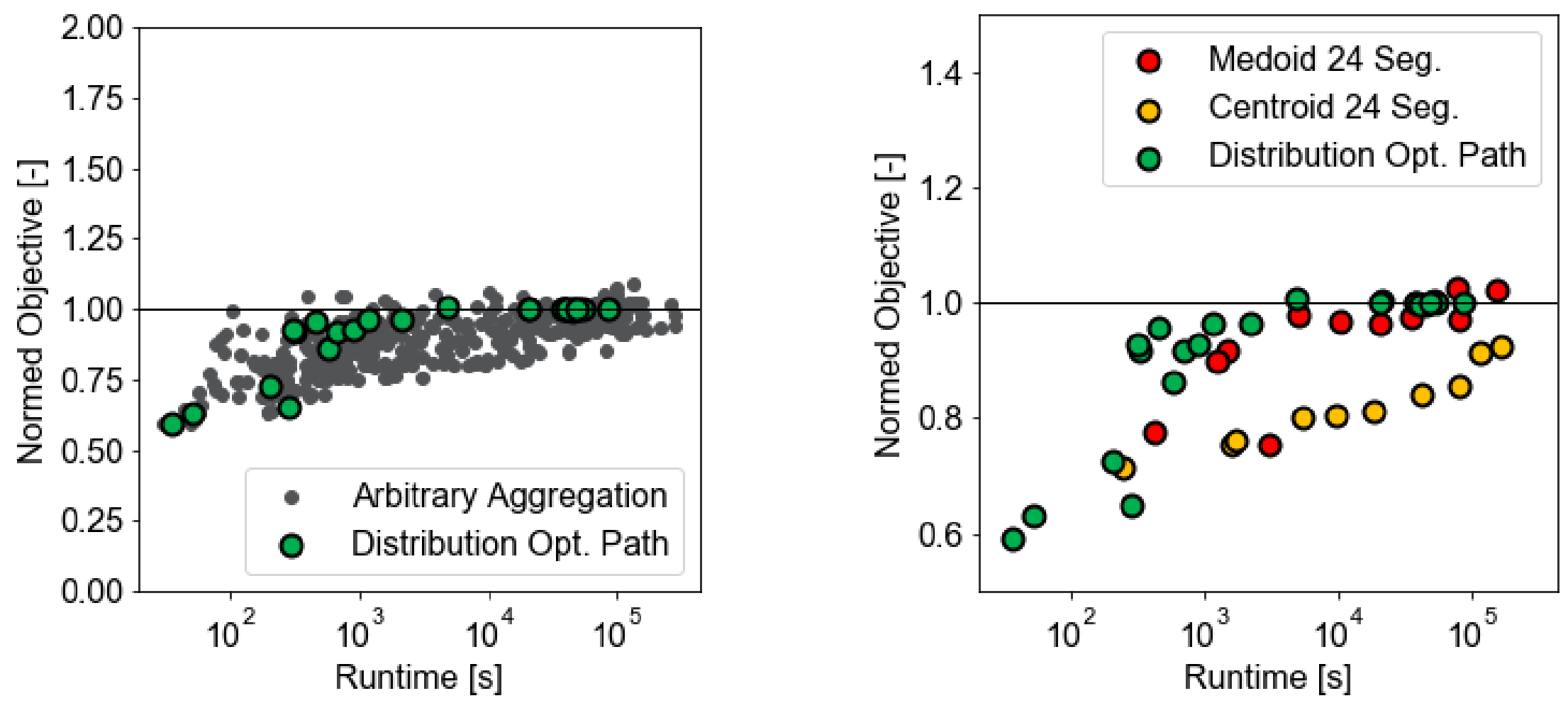

*Figure 22. The identified optimal aggregation method among all tested aggregation approaches (left) and a direct comparison to state-of-the-art aggregation methods as used in FINE (right).*

## 4.3. Additional observations

Apart from the analyses that could be individually conducted for each model, some observations arise from the comparison of both models and the time series employed.

First, the European system model differs considerably from the self-sufficient building one with respect to the number of time series, the number of regions, and the fact that numerous aperiodic wind time series are considered. Therefore, it can be used to study the effect of different temporal aggregation procedures on the technology choice. As a major part of the time series is given by the capacity factor time series for wind turbines and photovoltaic panels, these technologies are of particular interest.

Therefore, Figure 23 depicts the energy supply and storage capacities depending on a varying number of typical days and segments using a representation by centroids. Considering the configuration in the first row with 96 typical days and 12 segments per day a highly resolved model, a decrease to two segments per typical day, as shown in the second row, leads to a significant oversizing of the photovoltaic capacities. The reason for this is the disproportionate underestimation of solar intermittency in this case, which makes photovoltaic panels more economically attractive. In contrast, twelve typical days with twelve segments per day leads to a predominant underestimation of the wind time series' variance. Accordingly, the capacities for wind turbines are overestimated, as can be seen in the third row of Figure 23. Apart from that, a small number of typical days leads to an underestimation of seasonality effects. For that reason, the capacities for seasonal hydrogen storage are underestimated in the corresponding right column.

Accordingly, temporal aggregation and resolution in general have a significant and heterogeneous impact on different technologies, which must be taken into account in the modeling process.

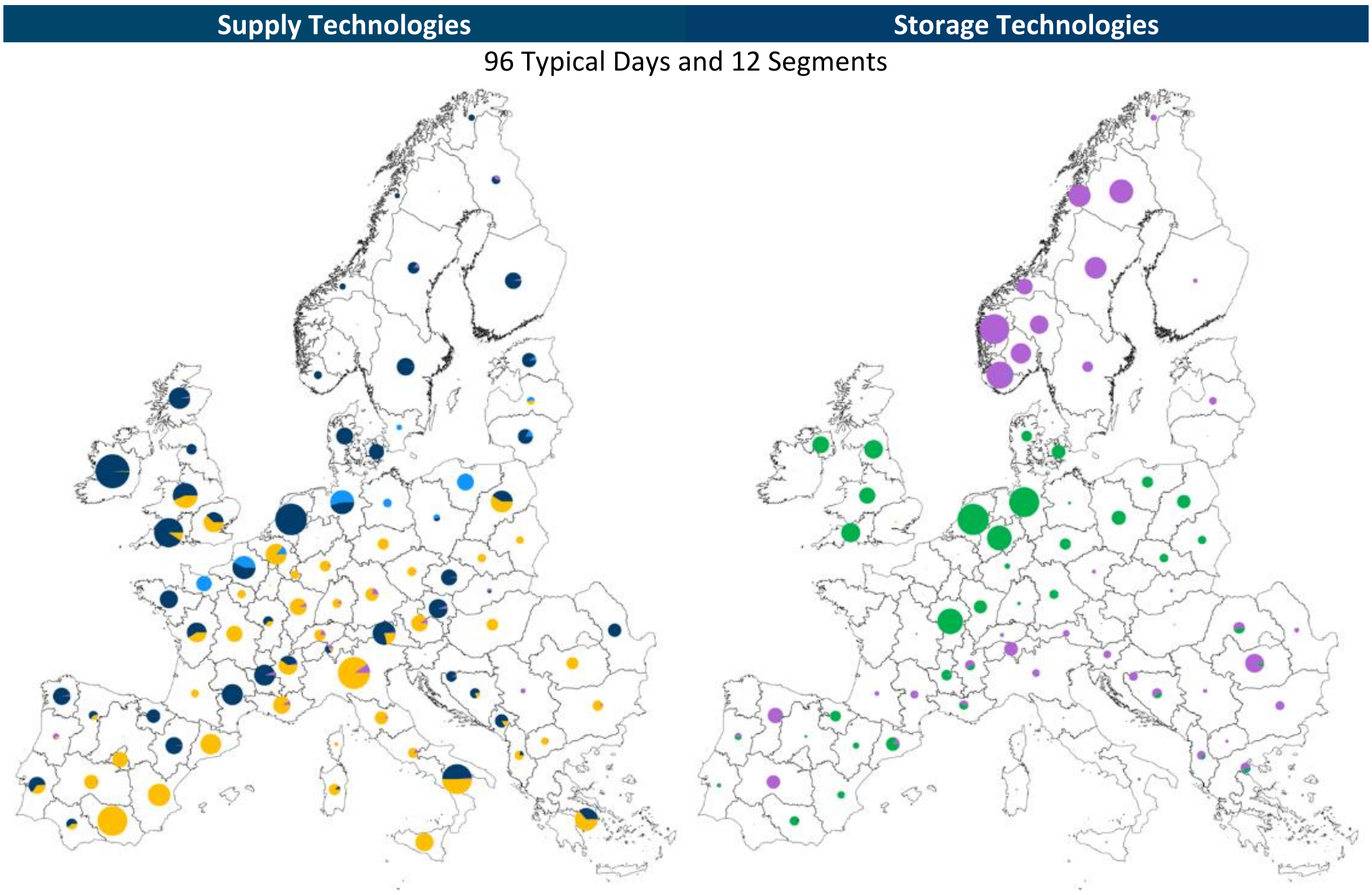

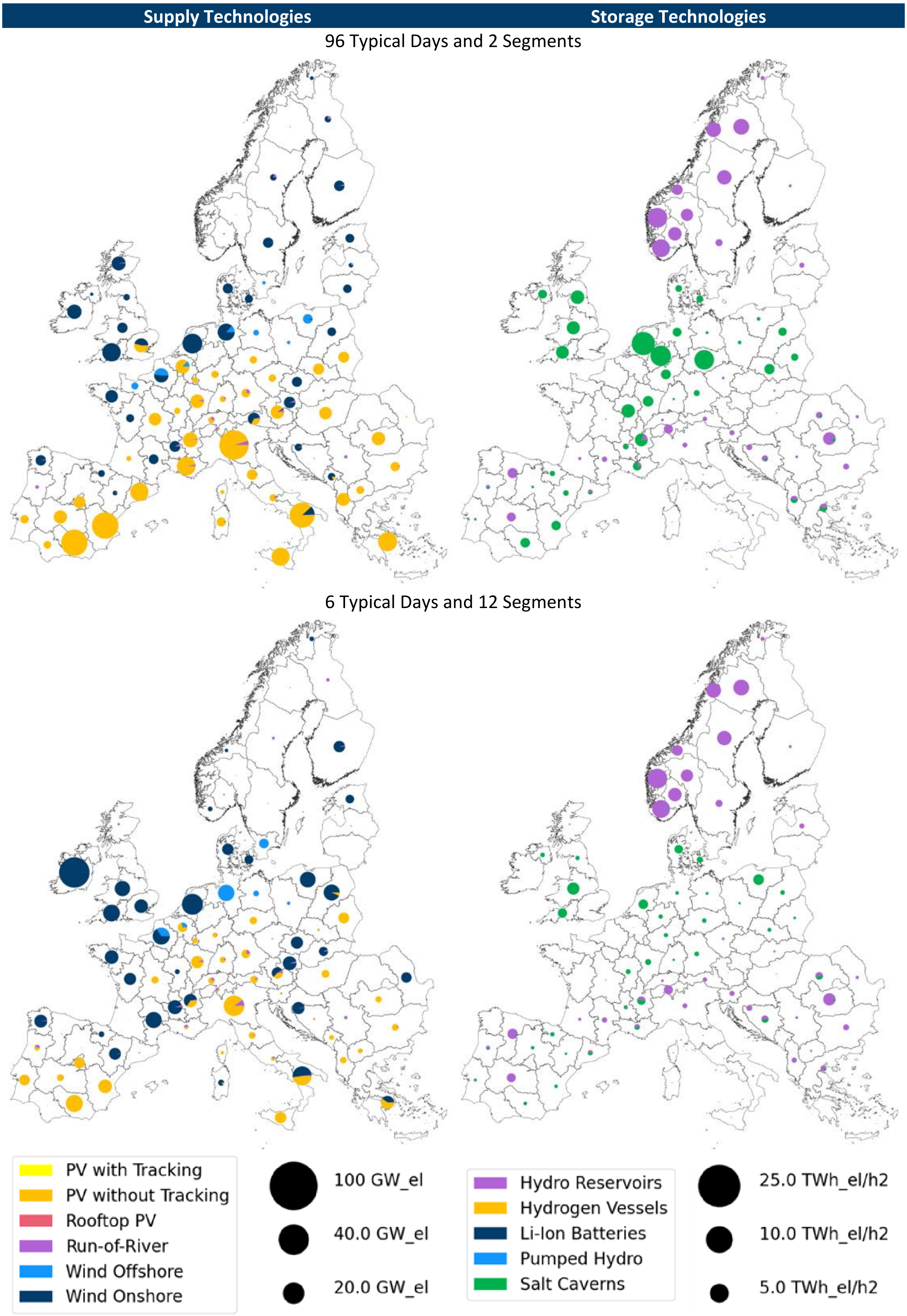

*Figure 23. Energy supply and storage capacities depending on a varying number of typical days and segments using a representation by centroids.*

Apart from that, both models vary with respect to the overall number of considered time series. Accordingly, it can be expected that the number of "redundant" time steps, which can be aggregated without a major loss of accuracy, may vary. For that, Figure 24 depicts the convergence behavior of the two models for an increasing number of typical days, with hourly time steps represented by centroids. The x-axis represents the total number of time steps of the respective aggregation configuration, with 8760 (365 typical days) serving as the reference case.

Due to the logarithmic scale of the x-axis, the marginal gain in accuracy decreases with an increasing number of typical days. However, the normed objective function value of the reference case forms an asymptote in the case of the self-sufficient building, which is not the case for the European model. For more than 1000 time steps, i.e., about 40 typical days, there is little change in the objective function value, which means that about 90% of the original time series values are indeed redundant. In contrast to that, the results for 1000 time steps are still far from those of the reference case for the European model, and it can be expected that the objective function value would further increase if more than 365 days were considered as original inputs.

The reason for this observation can be found in the multidimensionality of the respective datasets: In the case of the European model, much more low-correlated time series are simultaneously clustered and aggregated. Therefore, the resemblance of different time steps to each other is much lower, as a close similarity of two separate time steps becomes less likely as more low-correlated time series are taken into account. From a practical point of view, this also indicates that the number of potential extreme situations increases with the overall system size.

Therefore, the quality of temporal aggregation also depends on the redundancy, and ultimately the number and cross-correlation of the considered time series. This must be taken into account when choosing an appropriate method for reducing a model's mathematical complexity.

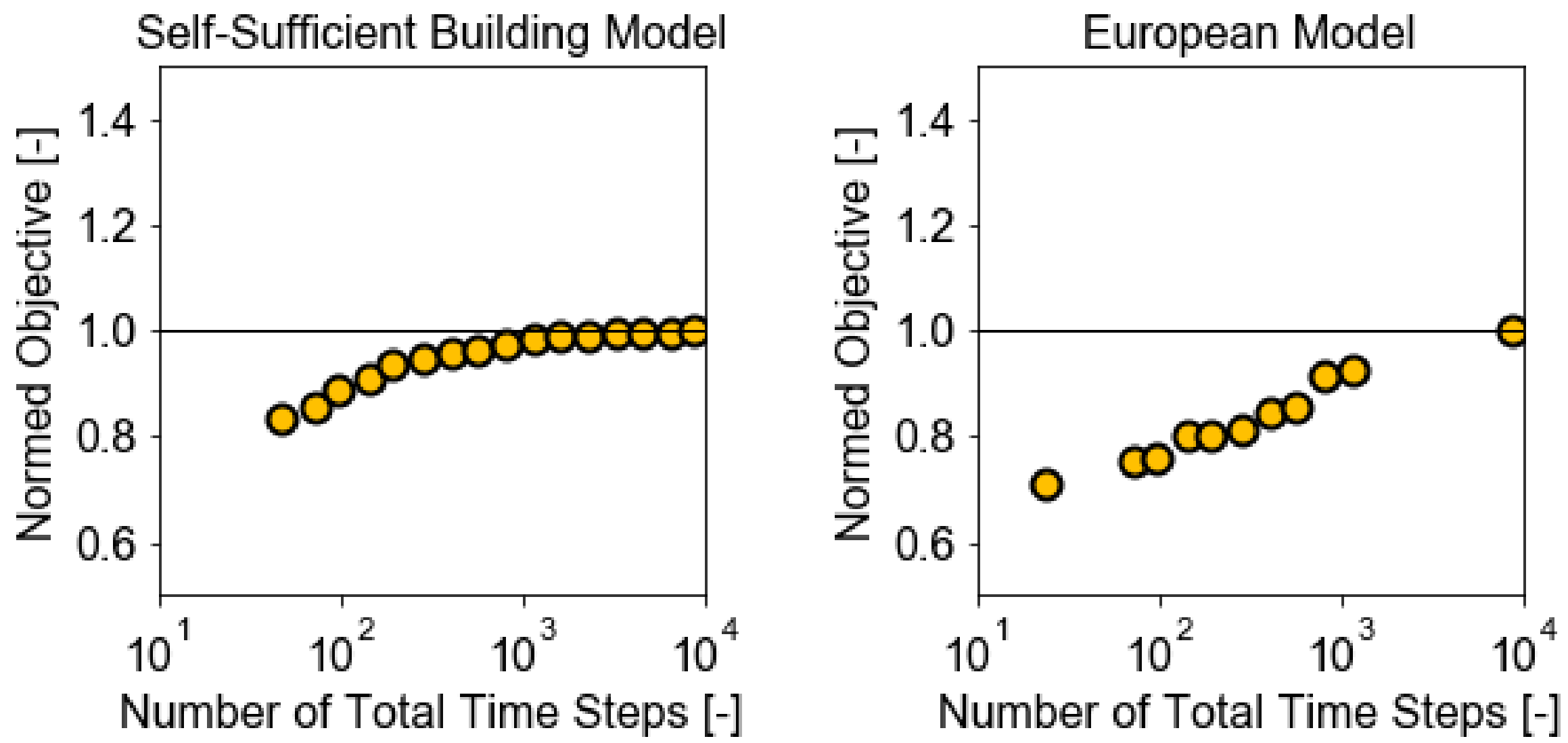

*Figure 24. Convergence behavior of the two different aggregated energy systems' optimal objective values to the reference value if centroid-based typical days with an hourly resolution (no segmentation) are used.*

## 5. Conclusions

In this study, a broad sensitivity analysis for different temporal aggregation techniques was applied to two different models and statistically-analyzed with respect to the desirable properties of an aggregation procedure, which led to small deviations in the respective energy system models. Concisely, two hypotheses were deduced stating that an imbalanced ratio between the number of typical time steps and the inner-daily temporal resolution may systematically bias the appropriate modeling of certain technologies, and that the preservation of the original time series' statistical momentums, and therefore its individual value distribution, is of the utmost importance for a good temporal aggregation. The individual research questions presented in the introduction of this study can be answered as follows:

**Research question 1:** Is there an optimal tradeoff between the duration of time steps and the number of typical periods, and are certain technologies affected systematically by either of these temporal aggregation approaches?

Our results imply that *there is* a systematic bias introduced by too low a number of inner-daily segments or typical days. We found that too low an inner-daily temporal resolution disproportionally underestimates the daily variance of photovoltaic time series and makes them more economically viable due to this underestimated intermittency. In contrast, too small a number of typical days more strongly affects the aperiodic behavior of wind profiles, with the same effect of overestimating wind capacities if only a few typical days are considered. We proposed an algorithm that automatically differentiates between the different time series and searches for a good tradeoff between the number of typical days and inner-daily segments for a given maximum number of total time steps.

**Research question 2:** Can the aggregation procedure be improved by explicitly taking statistical momentums, i.e., the original time series' value distributions, into account?

We developed an algorithm that is capable of aggregating time series in such a way that the original time series' value distribution is very closely approximated. The algorithm significantly outperformed a representation by centroids and medoids for both models. Therefore, it can be stated that not only the mean of a time series but also other statistical momentums of higher order are of importance for achieving an appropriate aggregation.

Both of the approaches proposed for improving on current temporal aggregation techniques led to a Pareto-optimal tradeoff between the calculation runtime and the energy system models' aggregation-induced objective function deviation. This means that for arbitrarily strongly aggregated, and therefore accelerated models, the results using the proposed methods were strictly more accurate than any comparably strongly aggregated model using traditional methods. On average, the speedup to traditionally-employed aggregation approaches was about one order of magnitude given the same accuracy level. Future research could focus on an improvement in the algorithms to account for the importance of certain time series for the energy system model to be aggregated when determining a ratio between the number of segments and typical days, as well as the development of a representation method that also preserves the original set of time series' covariance.

The methods are open source and available at: https://github.com/FZJ-IEK3-VSA/tsam.

## Acknowledgements

This work was supported by the Federal Ministry for Economic Affairs and Energy of Germany as part of the METIS project (project number: 03ET4064A).

# A. Appendix

## A.1. Techno-economic data of the self-sufficient building

*Table 4. Cost parameters of the self-sufficient building model.*

| Components | Capex | | | | Opex | | Lifetime | | Source |
|---|---|---|---|---|---|---|---|---|---|
| | Fixed | | Capacity-Specific | | Fixed+Capacity-Specific | | | | |
| Photovoltaic Ground | — | — | 4000.00 | €/kW$_p$ | 1.00 | % Inv./a | 20 | a | own assumption |
| Photovoltaic Rooftop | — | — | 769.00 | €/kW$_p$ | 1.00 | % Inv./a | 20 | a | [102] |
| Inverter | — | — | 75.00 | €/kW$_p$ | — | — | 20 | a | [103] |
| Battery | — | — | 301.00 | €/kWh$_p$ | — | — | 15 | a | [102] |
| Reversible Solid Oxide Cell | 5,000.00 | € | 2,400.00 | €/kW$_{el}$ | 1.00 | % Inv./a | 15 | a | [104] |
| Heat Pump | 4,230.00 | € | 504.90 | €/kW$_{th}$ | 1.50 | % Inv./a | 20 | a | [105] |
| Thermal Storage | — | — | 90.00 | €/kWh$_{th}$ | 0.01 | % Inv./a | 25 | a | [106] |
| E-Heater & E-Boiler | — | — | 60.00 | €/kW$_{th}$ | 2.00 | % Inv./a | 30 | a | [106] |
| Tank | — | — | 0.79 | €/kWh$_{H2}$ | — | — | 25 | a | [107] |
| Dibenzyltoluene | — | — | 1.25 | €/kWh$_{H2}$ | — | — | 25 | a | [108, 109] |
| Hydrogen Vessels | — | — | 15.00 | €/kWh$_{H2}$ | — | — | 25 | a | [110] |
| Hydrogenizer | 2,123.30 | € | 761.10 | €/kW$_{H2}$ | 1.00 | % Inv./a | 20 | a | [108] |
| Dehydrogenizer | 1,140.00 | € | 408.60 | €/kW$_{H2}$ | 1.00 | % Inv./a | 20 | a | [108] |
| Low Pressure Compressor | — | — | 1716.71 | €/kW$_p$ | 1.00 | % Inv./a | 25 | a | [111] |
| High Pressure Compressor | 560.00 | € | 1329.80 | €/kW$_p$ | 1.00 | % Inv./a | 25 | a | [111] |
| Heat-Exchangers 1 and 2 | — | — | 1.00 | €/kW$_{th}$ | 1.00 | % Inv./a | — | a | own assumption |
| Expanders 1 and 2 | — | — | 1.00 | €/kW$_{th}$ | 1.00 | % Inv./a | 25 | a | own assumption |